\def\AA         {{\bf A}}
\def\ZZ         {{\bf Z}}

\def\CC         {{\bf C}}
\def\QQ         {{\bf Q}}
\def\PP         {{\bf P}}
\def\ii         {{\rm i}}
\def\ee         {{\rm e}}

\def\vac        {\hspace{-4pt}>}

\def\Fock       {{\rm Fock}}

\def\MSV        {chiral de Rham complex}
\def\msv        {{\cal MSV}}

\def\dim        {{\rm dim}}
\def\sin        {{\rm sin}}
\def\Sym        {{\rm Sym}}
\def\sinh       {{\rm sinh}}
\def\deg        {{\rm deg}}
\def\Ell        {{\cal ELL}}

\def\Box        {{\bullet}}

\documentstyle[twoside,12pt]{article}
\newtheorem{prop}{Proposition}[section]
\newtheorem{dfn}[prop]{Definition}   
\newtheorem{theo}[prop]{Theorem}

\newtheorem{rem}[prop]{Remark} 
 
\newtheorem{lem}[prop]{Lemma}

\title{Elliptic Genera and Applications to Mirror Symmetry}

\author{
Lev A. Borisov
\\
\small Department of Mathematics,  Columbia University, New York, NY
10027\\
\small e-mail: lborisov@math.columbia.edu\\
Anatoly Libgober\\
\small Department of Mathematics, University of Illinois, Chicago, IL
60607\\
\small e-mail:libgober@math.uic.edu}

\begin{document}

\date{} 

\maketitle

\begin{abstract}
{The paper contains a proof that elliptic genus of a Calabi-Yau 
manifold is a Jacobi form, finds in which dimensions the elliptic 
genus is determined by the Hodge numbers and shows that
elliptic genera of a Calabi-Yau hypersurface in a toric variety
and its mirror coincide up to sign. 
The proof of the mirror property is based on the extension of 
elliptic genus to Calabi-Yau hypersurfaces in toric varieties 
with Gorenstein singularities. 
}
\end{abstract}

\section{Introduction}
One of the motivations for this paper was an attempt to 
understand those invariants of Calabi-Yau manifolds 
which value on the mirror $X^*$ is determined by the 
value on original manifold $X$.
Examples of such invariants that attracted the most attention 
are topological Euler characteristic,
Hodge numbers and various $d$-point functions. 
In particular the relation between the 
different kinds of $d$-point functions 
yield the famous predictions for enumerative geometry (cf. \cite{CoxKatz}).

The invariant considered in this paper is elliptic genus for 
Calabi-Yau manifolds. Elliptic genus of oriented differentiable 
manifolds first appeared in the works of Landweber-Stong and Ochanine
(cf. \cite{Landweber.Stong} and references there) 
as part of attempts to find genera satisfying certain topological 
conditions and also as a mean for constructing new generalized 
cohomology theories (elliptic cohomology). Elliptic genus 
defined in such way is a certain homomorphism from the
ring of oriented cobordisms $\Omega^*_{SO}$ 
into the ring of modular forms 
for the group $\Gamma_0(2)=\{ \pmatrix {a & b \cr c &d \cr}
\vert c \equiv 0\,
{\rm mod} 2 \}$. The value of 
elliptic genus at the cusps 
of $\Gamma_0(2)$ is equal to $\hat A$-genus and the signature.
At the same time Witten proposed a description of elliptic genus 
as an index of certain Dirac-like operator on the loop space
or more generally as the trace of certain operator associated 
with a conformal field theory (cf. \cite{Landweber.Stong}). 
Also at the same time Witten proposed a generalization of elliptic genus
for almost complex manifolds with first Chern class divisible by $N$ 
(cf. \cite{Landweber.Stong}). A construction of such generalization was
also given  independently by F.Hirzebruch (cf. \cite{Hirz2}). These genera
are modular forms for $\Gamma_1(N)$. Motivated by attempts to extend  
rigidity property of elliptic genera, Krichever (cf. \cite{Krichever})
introduced an extension of Witten-Hirzebruch genus. Values of
Krichever's elliptic genus are certain functions of several
variables. 
 
Also, in the physics literature there was proposed a version of elliptic
genus such that its values on Calabi-Yau manifolds are (weak) Jacobi
forms, i.e. certain functions on $H \times {\bf C}$. Several interesting
properties were observed (cf. \cite{DMVV,EOTY,KYY}). One of
the most striking is the calculation of the elliptic genus 
on symmetric powers of a manifold $M$ is terms of the elliptic genus of
$M$ itself (cf. \cite{DMVV}). Also several comments were made about
possible relationship with mirror symmetry.

The purpose of this paper is to address from the mathematical point of
view the properties of  elliptic genera of Calabi-Yau manifolds and
to study elliptic genera of toric varieties. We start by reviewing in
Section 2 some aspects of previous work by Landweber-Stong, Ochanine,
Witten, Hirzebruch and Krichever mentioned above. 
Then in Section 3 we introduce elliptic genus for Calabi-Yau manifolds,
following suggestions from physics literature as the holomorphic Euler
characteristic of the bundle:
 $$y^{-{{dim M} \over 2}} \otimes_{n \ge 1} (\Lambda_{-yq^{n-1}}
\bar T_M \otimes \Lambda_{-y^{-1}q^n} T_M \otimes S_{q^n}\bar T_M
\otimes S_{q^n}T_M)$$  
While this Euler characteristic can be calculated for 
any manifold and yields a holomorphic function of $H \times {\bf C}$, 
we show that for Calabi-Yau manifolds this Euler characteristic
is a Jacobi form of weight $0$ and index $d \over 2$ where $d$ is the
dimension of the manifold. This is done by using an
expression for the characteristic series for such 
genus via theta functions. Theta functions 
did appear already in connection with elliptic genus
but for the proofs of rigidity properties (cf. \cite{Liu}) .
We also show how elliptic genera of differentiable manifolds
and almost complex manifolds mentioned above are related to the
elliptic genus considered in present paper in the case of Calabi-Yau
manifolds.

The crucial question is whether elliptic genus of Calabi-Yau 
manifold can be expressed in terms of Hodge numbers, i.e. 
whether the behavior of the elliptic genus in mirror 
correspondence can be deduced from known behavior of the 
Hodge numbers. We show in this paper that the answer is 
"no'' provided the dimension of the manifold is at least 12 (15 for
odd dimension). The Jacobi property of elliptic genus is the key issue in
the argument. It allows us to calculate dimensions of 
spaces of functions on $H \times {\bf C}$ which are the 
elliptic genera of Calabi-Yau manifolds either by 
interpreting Jacobi forms as sections of certain bundles on 
the compactification of the quotient of $H \times {\bf C}$
by the Jacobi group or by  using calculations of the space of 
weak Jacobi forms due to Eichler-Zagier (cf. \cite{EZ}). A consequence of
this is it that the space of functions generated by elliptic genera of
Calabi-Yau manifolds is a free algebra on three generators of degrees 1, 2
and 3. This also implies that the space of elliptic genera of manifolds of
dimension greater than 13 (or is equal to 12) is too big and elliptic 
genera cannot just depend on Hodge numbers. All of this is accomplished in
Section 4.

Next we consider the elliptic genera for toric varieties.
Use of the torus action of the bundle in the definition 
of the elliptic genus yields a series representing the 
elliptic genus in terms of the defining fan. One can restrict 
this elliptic genus to $y=-1$ which yields a formula for the 
elliptic genus studied by Landweber-Stong-Ochanine and Witten.
Comparison of these series with standard expressions in 
terms of Eisenstein series yield very interesting identities.
For example for $M={\bf P}^2$ (after some easy modifications) we obtain
the following identity:
 $$\sum_{m \ge 1,n \ge 1} {{q^{m+n}} \over {(1+q^n)(1+q^m)(1+q^{m+n})}}
=\sum_{r \ge 1}q^{2r}\sum_{k \vert r} k$$
One can prove this identity by elementary means (cf. Section 5) but huge
class of identities corresponding to other toric manifolds is somewhat
mysterious.

Finally, in the last two sections the relationship 
  $$Ell(X;y,q)=(-1)^d Ell (X^*;y,q) \eqno (*)$$ 
between elliptic genera of a $d$-dimensional Calabi-Yau hypersurface $X$
in a Fano toric variety and its mirror $X^*$ is derived. The proof relies
heavily on the work \cite{Bvertex} by the first author. There are several
ingredients in it which we hope have independent interest. First, the
starting point is an interpretation of elliptic genus as certain trace
which is a reminiscence of original Witten's definition but the
(super)trace here is calculated on the cohomology of the chiral de Rham
complex studied in \cite{MSV} and the work \cite{Bvertex} by the first
author. The material from \cite{Bvertex} needed for the proofs here is
reviewed and used in Section 6. Second, since the chiral de Rham
complex was defined in \cite{Bvertex} for Gorenstein toric varieties and
Calabi-Yau hypersurfaces in Gorenstein toric Fano varieties this trace
formula allows to define elliptic genus for such singular varieties as
well. Based on results from \cite{Bvertex} transformation law for the
elliptic genus under mirror correspondence is proven. Then in the Section
7 it is  shown that the elliptic genus of a Calabi-Yau hypersurface with
Gorenstein singularities is a weak Jacobi form. This, together with the
transformation law yields the relation (*).

\section{A short review of elliptic genus}
A genus (resp. complex genus) with values in a $\bf Q$-algebra 
 $R$ with unit 
is a ring homomorphism from the oriented cobordism ring $\Omega^{S0}_*$
(resp. complex cobordism ring $\Omega^{U}_*$) to $R$. 
According to \cite{Hirz1} such homomorphisms are in 
one to one correspondence with the formal power series $Q(x)$ with 
coefficients in $R$ satisfying $Q(0)=1$. The genus $\Psi_Q(x)$ 
corresponding to a series $Q(x)$ can be described as follows. 
Let $c(X)=\prod (1+x_i)$ be a formal factorization of the total 
Pontryagin (in complex case Chern class) of a manifold $X$.  
Then the genus $\Psi_q(X)$ is: 
$$
\prod Q(x_i) [X] 
$$
where $[X]$ is the fundamental class of $X$ and $\prod Q(x_i)$ is 
written as a polynomial in symmetric functions in $x_i$ i.e. Pontryagin
(resp. Chern in complex case) classes of $X$. 
If $Q(0) \ne 1$ but $Q(0) \ne 0$ 
then  the above formula still produces a $R$-valued invariant of 
the manifold. This ``non-normalized'' genus $\Psi_{Q}(X)$ 
is related to the genus $\Psi_{\tilde Q}(X)$,  corresponding to
the series  $\tilde Q(x)={{Q(x)} \over {Q(0)}}$
as follows: $\Psi_Q(X)=\Psi_{\tilde Q}(X) 
\cdot Q(0)^{{1 \over 4}dim_{\bf R}(X)}$ (exponent is $dim_{\bf C}X$ 
in complex case). 

Elliptic genus of an oriented manifold $X$ can be defined as 
${\bf Q}[[q]]$-valued genus corresponding to the series
(cf. volume \cite{Landweber.Stong} and references there):
$$
Q(x)={{x/2}\over {\sinh( x/2)}} \prod_{n=1}^{\infty}
  [{{(1-q^n)^2} \over {(1-q^n\ee^x)(1-q^n\ee^{-x})}}]^{(-1)^n} 
$$
It can be described also as  
     $$\hat A(X) ch \{ {{R(T)} \over {R(1)^{dim X}}} \}[X]  $$
Here $T$ is the tangent bundle, $\hat A$ is the ($\bf Q$-valued) 
genus corresponding to the series ${{x/2}\over {\sinh x/2}}$ and
   $$ R(T)=\otimes_{l>0, l \equiv 0 (2)} S_{q^l}(T) 
   \otimes_{l>0, l \equiv 1 (2)} \Lambda_{q^l}(T) $$
where $$S_q(V)=\Sigma S^n(V) q^n,  \Lambda_q(V)=\Sigma \Lambda^n(V)q^n $$ 
are generating series for symmetric and exteriors powers of a bundle 
$V$. 

Elliptic genus of an oriented manifold is a modular form 
on $\Gamma_0(2)$ with rational coefficients if $q=\ee^{2 \pi \ii  \tau}$
where $\tau$ is in the upper half plane $H$.

Hirzebruch (\cite{Hirz2}) and Witten (\cite{Landweber.Stong})
defined genera of complex manifolds which are modular forms
on:
 $$\Gamma_1(N)=\{ \pmatrix{a \ b \cr c \ d\cr}
\in SL(2,{\bf Z})
\vert c \equiv 0 (N), a \equiv d \equiv 1 (N) \} $$ 
provided the first Chern class of the manifold satisfies 
$c_1 \equiv 0 (N)$. 

In terms of characteristic series these genera 
can be defined as follows. For  $\tau \in H$ 
$$\alpha=2 \pi \ii ({k \over N} \tau+ {l \over N}) \ne 0$$
 and 
\begin{equation}
\Phi(x, \tau)=(1-\ee^{-x}) \prod_{n=1}^{\infty}
  {{(1-q^n\ee^x)(1-q^n\ee^{-x})} \over {(1-q^n)^2}} 
\label{hirzebruch}
\end{equation}
the characteristic series is:
 $$Q_{HW}(x,\tau)=x \ee^{-{k \over N}x} {{\Phi(x-\alpha)} \over {\Phi(x)
     \Phi(-\alpha)}}$$
 
Krichever (\cite{Krichever}) considered a genus with characteristic
series 
$$
Q_K(x,z,\omega_1,\omega_2,\kappa)=x\ee^{-\kappa x}
{{\sigma_{\omega_1,\omega_2}(x-z)} \over 
{\sigma_{\omega_1,\omega_2}(x)
    \sigma_{\omega_1,\omega_2}(-z)}}\ee^{\zeta_{\omega_1,\omega_2}(z)x}
$$
where $z, \kappa \in {\bf C}^*$, 
$\sigma_{\omega_1,\omega_2}(z)$ and $\zeta_{\omega_1,\omega_2}(z)$ 
are Weierstrass 
$\sigma$ and $\zeta$ functions corresponding to the same lattice 
in $\bf C$ (cf. \cite{Chandra}). If the lattice is $L=2 \pi \ii \tau{\bf
Z} +2 \pi \ii {\bf Z}$ then 
the series $Q_K$ specializes into $Q_{HW}$ for  $z=\alpha$ and 
$\kappa=-{{2k} \over N} \zeta({\pi \ii \tau})-{{2l} \over N} \zeta (\pi
\ii)+\zeta (z)$. 

Krichever proved rigidity theorem for such genus i.e. showed that 
for SU-manifolds with $S^1$-action the corresponding equivariant
genus is a multiple of a trivial character  generalizing  
similar results for orientable and complex manifolds 
(cf. \cite{Taubes,Hirz2}). 

In recent preprint (\cite{Totaro}) Burt Totaro identified the image of the
universal  genus corresponding to the series $Q_K$ as the quotient of the 
SU-cobordism ring by the equivalence relation generated by (classical)
flops.

In \cite{EOTY},\cite{DMVV}  the authors considered a genus for (almost) 
complex manifolds for which $c_1=0$ or equivalently the structure group 
of the
 tangent bundle can be reduced to the group $SU$. It can be
defined either as a (super)trace of a certain operator 
(cf. Definition 7.1 below) or as
\begin{equation} 
\chi(M,q,y)=\int_Mch(\Ell_{q,y})td(M)
\label{elliptic1} 
\end{equation} where

\begin{equation} 
\Ell_{q,y}=y^{-{\dim M \over 2}} \otimes_{n \ge 1} (\Lambda_{-yq^{n-1}}
\bar
T_M \otimes \Lambda_{-y^{-1}q^n} T_{M} \otimes S_{q^n}\bar T_M \otimes
S_{q^n} T_M) \label{elliptic2} \end{equation}

\bigskip 
In the next section we show that for Calabi-Yau 
manifolds the above 
elliptic genus is a weak Jacobi form. 
While this result is stated in \cite{KYY},
there seems to be no proof of it anywhere in the literature. 
We also spell out the relationship between the genus 
(\ref{elliptic1}),(\ref{elliptic2}) and the genera corresponding to 
the series $Q_{HW}$ and $Q_K$. As a consequence we see that 
for Calabi-Yau manifolds  the above elliptic genus up to simple factors
coincides with the genus corresponding to $Q_K$.

\section{Elliptic genera of Calabi-Yau varieties as weak Jacobi forms}

In this section we assume that $M$ is a complex compact manifold 
of dimension $d$.

Let $$\theta(z,\tau)=q^{1 \over 8}  (2 \sin \pi z)
\prod_{l=1}^{l=\infty}(1-q^l)
 \prod_{l=1}^{l=\infty}(1-q^l \ee^{2 \pi \ii z})(1-q^l \ee^{-2 \pi \ii
z})$$
where $q=\ee^{2 \pi \ii  \tau}$ is the Jacobi's theta function
(cf. \cite{Chandra}). As a theta function with characteristic this
 is $\theta_{1,1}(z,\tau)$ (cf. \cite{Mumford}).
We have 
$$\theta({z \over \tau},-{1 \over \tau})=
 -\ii \sqrt { \tau \over \ii }
\ee^{{\pi \ii z^2} \over {\tau}} \theta(z,\tau)$$
(cf. \cite{Chandra}).
\bigskip

\begin{prop}
{\rm 
If $c(T_M)=\prod(1+x_i)$ then
$$\chi(M,q,y)=\int_Mch(\Ell_{q,y})td(M)$$ where
$$\Ell_{q,y}=y^{-{d \over 2}} \otimes_{n \ge 1} \Lambda_{-yq^{n-1}}\bar 
T_M \otimes \Lambda_{-y^{-1}q^n} T_{M} \otimes S_{q^n}\bar T_M \otimes
S_{q^n} T_M)$$ is equal to the integral over $M$ of
the degree $d$ term in the expansion of:
$$
E(z,\tau, x_1,...,x_r)=
   \prod_i x_i {{\theta ({{x_i} \over {2 \pi \ii}}-z,\tau)} \over
{\theta ({{x_i} \over {2 \pi \ii }}, \tau)}} 
$$
where $y=\ee^{2 \pi \ii z}$.
}
\label{claim1}
\end{prop}

{\em Proof.}
Indeed, $\chi$ is the genus corresponding to the series:
$$y^{-{1 \over 2}}\prod_{n=1}^{n=\infty}
{{(1-yq^{n-1}\ee^{-x})(1-y^{-1}q^n\ee^{x})} \over
{(1-q^n\ee^{-x}) (1-q^n\ee^{x})}} \cdot {x \over {1-\ee^{-x}}}=$$
$$y^{-{1 \over 2}}\prod_{n=1}^{n=\infty}
{{(1-yq^{n}\ee^{-x})(1-y^{-1}q^n\ee^{x})(1-y\ee^{-x})} \over
{(1-q^n\ee^{-x}) (1-q^n\ee^{x})}} \cdot {x \over {1-\ee^{-x}}}.$$
Hence for $M$ with factored Chern class we have the following  
expression for generating series for $\chi$
$$\ee^{-\pi \ii d z} \prod_{n=1}^{n=\infty}
{{(1-\ee^{2 \pi \ii (z-{{x_i} \over {2 \pi \ii }})}q^n)
(1-\ee^{-2 \pi \ii (z-{{x_i} \over {2 \pi \ii }})}q^n)}
\over {(1-\ee^{2\pi \ii {{x_i } \over {2 \pi \ii}}}q^n)
(1-\ee^{-2 \pi \ii  {{x_i} \over {2 \pi \ii}}}q^n)}} \cdot  
 {{(1-\ee^{2 \pi \ii z-{x_i} }) x_i} \over
 {(1-\ee^{-x_i})}}=$$
$$\ee^{-\pi \ii d z}\prod_i {{\theta(z-{{x_i} \over {2 \pi \ii }})
\sin (-\pi {{x_i} \over {2 \pi \ii }})} \over {\theta( -{{x_i} \over {2
\pi
i
}}) \sin \pi (z-{{x_i} \over {2 \pi \ii }})}} \cdot
{{(1-\ee^{2 \pi \ii z- x_i}) x_i} \over {(1-\ee^{-x_i})}}=$$
$$\ee^{-\pi \ii d z} \prod_i {\theta(z-{{x_i} \over {2 \pi \ii }}) \over
{\theta(-{{x_i} \over {2 \pi \ii}})}}
\prod_i {(\ee^{{-x_i} \over 2}-\ee^{{{x_i} \over 2}}) \over
 (\ee^{\pi \ii z -{{x_i} \over 2 }}-\ee^{-\pi \ii z+{{x_i} \over 2}}) }
\cdot{{(1-\ee^{2 \pi \ii z-x_i}) x_i } \over {(1-\ee^{-x_i})}}=
$$
$$
\prod_i {{x_i \theta(z-{{x_i} \over {2 \pi \ii }})}
\over {\theta (-{{x_i} \over {2 \pi \ii }})}}=
\prod_i {{x_i \theta({{x_i} \over {2 \pi \ii }}-z)}
\over {\theta ({{x_i} \over {2 \pi \ii }})}}.$$\\[-2em]   
{$~$}$\hfill\Box$

\begin{theo}
{\rm 
Function $\chi (M,z, \tau)$ is a (weak) Jacobi
of weight $0$ and index $d/2$.} 
\end{theo}

\bigskip 
Weak here implies that while it obeys the transformation laws of the
Jacobi forms, it does not satisfy regularity conditions at the cusps.
Rather, the only condition at the cusp is that $q$ appears with
nonnegative powers only. We refer to \cite{EZ} for precise definitions.
Also, when $d$ is odd, the definition of the Jacobi form must be modified
to allow a character (\cite{KYY}).
\bigskip 

{\em Proof.}
The condition at the cusp clearly holds because $\Ell_{q,y}$ has no
negative powers of $q$. To verify that $\chi (M,z, \tau)$ is a Jacobi form
of weight $0$ and index $d/2$ with character it is enough to check the
modular properties of $\chi(M,z,\tau)$ for generators of Jacobi group.
Indeed, if  $M$ is Calabi-Yau then we have:
\begin{equation} 
\chi(M,z,\tau+1)=\chi(M,z,\tau) 
\label{modular1} 
\end{equation}
\begin{equation}
\chi(M,{z \over \tau},-{1 \over \tau})=\ee^{{ \pi \ii d z^2} \over \tau}
\chi(M,z,\tau) \label{modular2} 
\end{equation}
\begin{equation} \chi(M,z+\tau,\tau)= (-1)^d
\ee^{- \pi \ii d (\tau + 2z)} \chi(M,z,\tau) \label{modular3}
\end{equation}
\begin{equation}
\chi(M,z+1,\tau)=(-1)^d \chi(M,z,\tau) \label{modular4} \end{equation}
(\ref{modular3}) and (\ref{modular4}) 
follows from the identities:   
 $$\theta (z+1,\tau)=-\theta(z,\tau), \ \ \ \theta(z+\tau,\tau )=-\ee^{-2
\pi \ii
z-\pi \ii \tau} \theta(z,\tau)$$ and (\ref{modular1}) is obvious.
Let $$\prod x_i  {{\theta({x_i \over {2 \pi \ii }}-z, \tau)}
\over {\theta ({{x_i} \over {2 \pi \ii }},\tau)}}= 
\sum_{\bf k} Q_{\bf k}(z,\tau) {\bf x}^{\bf k}$$
where $\bf x$ is a product of powers of $x_i$ and $\bf k$ is
multiindex. Hence, replacing $\tau \rightarrow {-{1 \over {\tau}}},
x_i \rightarrow {{x_i} \over \tau}$ and $z \rightarrow {z \over \tau}$
we obtain:
\begin{equation}\prod_i ({x_i \over {\tau}}) {{\theta (-{z \over
\tau}+{{x_i}
\over {2 \pi \ii \tau}},-{1 \over \tau})} \over
{\theta ({{x_i} \over {2 \pi \ii \tau}},-{1 \over \tau})}}=
\sum_{\bf k} Q_{\bf k}({z \over \tau},-{1 \over {\tau}})
 ({{\bf x} \over \tau})^{\bf k}  \label{modular5} \end{equation} 
Left hand side in this identity can be replaced by  
$$\prod_i({{x_i} \over {\tau}}){{\ee^{\pi \ii {{{(-z+{{x_i} \over
{2 \pi \ii}})^2}} \over {\tau}}} \theta (-z+{{x_i} \over {2 \pi \ii}},
\tau)}
\over {\ee^{\pi \ii {{({{x_i} \over {2 \pi \ii }})^2} \over \tau}}
\theta ({{x_i} \over {2 \pi \ii }})}}= $$
\begin{equation} ({1 \over \tau})^d \prod \ee^{-{z x_i} \over \tau}x_i
{{\ee^{{\pi \ii z^2} \over \tau}
\theta(-z+
{{x_i} \over {2 \pi \ii }})} \over {\theta ({{x_i} \over {2 \pi \ii }})}}
\label{modular6} \end{equation}
For degree $d$ terms in (\ref{modular5}) and (\ref{modular6}) we have
 $$Q_d({z \over \tau},-{1 \over \tau})=\ee^{{\pi \ii d z^2} \over \tau}  
 Q_d({z,\tau})$$ which yields (\ref{modular2}). $\hfill{\Box}$

\bigskip  
\bigskip

We shall conclude this section with a description of relationship 
between the genus considered in this section and genera corresponding to 
series $Q_{HW}$ and $Q_K$.

\bigskip 

\begin{prop}
{\rm Let $K(z,\omega_1,\omega_2,k)(X)$ be the 
genus of an almost complex manifold $X$ that 
corresponds to $Q_K$. Then for a smooth Calabi-Yau manifold $X$
of dimension $d$ 
$$ K(2\pi \ii z, 2 \pi \ii, 2 \pi \ii \tau,k)(X)=
Ell(z,\tau)(X) \cdot (-{{\theta^{\prime}(0,\tau)} 
\over {2\pi \ii\, \theta(z,\tau)}})^d
$$ 
where $\zeta_{\omega_1, \omega_2}$ is the Weierstrass
$\zeta$-function
corresponding to the lattice $\omega_1,\omega_2$.
}
\end{prop}

{\em Proof.}
Since (cf. \cite{Chandra} p.60)
$$\sigma_{\omega_1,\omega_2}(u)=\theta({u \over
{\omega_1}}, \tau)
\cdot {{\omega_1} \over {\theta}^{\prime}(0,\tau)} \cdot 
\ee^{\zeta_{\omega_1 ,\omega_2}(\omega_1/2) \cdot {{u^2} \over
{\omega_1}}}$$ 
where $\tau={{\omega_2} \over {\omega_1}}$
we have 
$$Q_K(x,2 \pi \ii z, 2 \pi \ii, 2 \pi \ii \tau, \kappa)=-x {{\theta ({x
\over
{2 \pi \ii }}-z)} \over {\theta ({x \over {2 \pi \ii }}) \theta(z)}} 
\cdot {{\theta^{\prime}(0,\tau)} \over {2 \pi \ii }}\cdot
\ee^{(-k+\zeta_{2\pi \ii ,2\pi\ii z }(2 \pi \ii z )
-2z\zeta_{2\pi \ii ,2\pi\ii z}(\pi\ii))x}$$
and the claim follows from $c_1(X)=0$ and Proposition \ref{claim1}.
\hfill$\Box$

\begin{prop} 
{\rm
Let $Ell(X,z,\tau)$ be elliptic genus and $H(X,\alpha,\beta,N,\tau)$ be
Hirzebruch elliptic genus of level $N$ ($\alpha,\beta \in {\bf Z}$ not
both divisible by $N$) for a $d$-dimensional complex Calabi-Yau manifold
$X$. Then for $\omega={{2 \pi \ii(\alpha \tau +\beta)} \over N}$ one has: 
$$\Psi_{HW} (M,\alpha,\beta,N,\tau)=
Ell(M,{\omega \over {2 \pi \ii }}, \tau) \cdot \theta(-{{\omega} 
\over {2 \pi \ii }})^{-d} \cdot \eta^{3d}(\tau)
$$
}
\end{prop}

{\em Proof.}
Using 
$$
\eta(\tau)=q^{1 \over 24} \prod_{l \ge 1} (1-q^l)
$$
the function  $\Phi$ given by (\ref{hirzebruch}) can be written as
$$\Phi(x, \tau)=(1-\ee^{-x})q^{1 \over 8}
{{\prod_{n=1}^{n=\infty}{(1-q^n\ee^x)(1-q^n\ee^{-x})(1-q^n)}} 
\over {\eta^3}}=$$
$$(1-\ee^{-x}) {{\theta({x \over {2 \pi \ii }}, \tau)}
\over {2 \sin (\pi {x \over {2 \pi \ii }}) \eta^3(\tau)}}=
\ee^{x \over 2} {\theta({x \over {2 \pi \ii }},\tau) \over {\eta^3(\tau)}}
$$
Therefore
$$Q_{HW}={{x  \Phi(x-\omega, \tau)} \over {\Phi(x,\tau)
\Phi(-\omega,\tau)}}\ee^{-{{\alpha x} \over N}}=
{{x\ee^{-{\alpha x} \over N} \theta({{x-\omega} \over {2 \pi \ii}}, \tau)
\eta^3(\tau)}
\over {\theta({x \over {2 \pi \ii }},\tau) \theta (-{\omega \over {2 \pi
i }},\tau)}} $$
and the claim follows from $c_1(X)=0$ and Proposition \ref{claim1}. 
\hfill$\Box$

\section{Spaces spanned by elliptic genera}
 \par {\it 0. Notations.} $\Gamma_{k,l}(n), k,l \in {\bf Q}, n \in {\bf
Z}$
is the semidirect product of $k{\bf Z} \oplus l{\bf Z}$ and $\Gamma (n)$
i.e. the collection of $(M,X)$ where $M \in \Gamma(n)$ and
$X=(x_1,x_2), x_1 \in k{\bf Z}, x_2 \in l{\bf Z}$. The product is
given by $(M,X)(M',X')=(MM',XM'+X')$. This product is the only
one for which $(\gamma_1 \cdot \gamma_2)(z,\tau)=\gamma_1(\gamma(z,\tau))$
with standard action of $\Gamma_{1,1}(1)$ on $H \times {\bf C})$
i.e. $(m,n,\gamma)(z,\tau)=({{z+m\tau+n} \over {c\tau+d}},\gamma(\tau))$.
\smallskip
\newline   ${\bf P}^2_n$ is the projective plane blown up at $n$ points.
\bigskip
\newline {\it 1. A compactification of} 
$H \times {\bf C}/\Gamma_{1,1}(1)$.

One can obtain a compactification for
$H \times {\bf C}/\Gamma_{1,1}(1)$ from compactification of quotients
corresponding to congruence subgroups which are particularly simple.

Let us consider the compactification of
$H \times {\bf C}/\Gamma_{1,1}(2)$ (cf. \cite{vanderGeer})
as the projective plane blown up at four base points of a generic 
pencil of quadrics. The compactification of
$H/\Gamma(2)$  is the base of the pencil of elliptic curves
(universal elliptic curve with level 2 structure) 
such that its quotient induced by involution $z \rightarrow -z$ 
is just mentioned pencil of quadrics. Three cusps of $\Gamma(2)$
correspond to three singular elements of this pencil.
The images of points of order 2 on the universal elliptic 
curve with level 2 structure in the pencil of quadrics 
form 4 sections of the pencil of quadrics which are 
four exceptional curves in ${\bf P}^2_4$. In particular each component of a 
singular member of the pencil intersects with two such sections, i.e. 
the components are in one to one correspondence with the pairs of 
distinct points of order 2.

Since $H \times {\bf C}/\Gamma_{1,1}(1)=(H \times {\bf
  C}/\Gamma_{1,1}(2))
/(\Gamma_{1,1}(1)/\Gamma_{1,1}(2))$ and 
$\Gamma_{1,1}(1)/\Gamma_{1,1}(2)=
SL_2({\bf F}_2)$ a compactification of
 $H \times {\bf C}/\Gamma_{1,1}(1)$ is given by 
the quotient    ${\bf P}^2_4/SL_2({\bf F}_2)$
(*). \footnote{(*) In fact it supports
a natural action of a bigger group $\Gamma_{{1 \over 2},{1 \over
2}}(1)/\Gamma_{1,1}(2)$ which is a semidirect product of ${\bf F}_2^2$ and
$SL_2({\bf F}_2)$ which is isomorphic to symmetric group $S_4$.The
action on ${\bf P}^2_4$ is induced by permutation of four blown up points.
The action of $S_3=SL_2({\bf F}_2)$ is the action of the subgroup of
$S_4$ fixing one of four sections of the pencil, i.e. the one
corresponding to the the zero of universal elliptic curve.} The
action of $\Gamma_{1,1}(1)/\Gamma_{1,1}(2)=SL_2({\bf Z})/\Gamma(2)=
SL_2({\bf F}_2)$ on singular fibres of the pencil of quadrics is specified
by the action of $SL_2({\bf Z})$ on points of order 2 on elliptic curves
and hence on the sections of the pencil of quadrics. One can equivariantly  
blow down three lines in singular fibers that do not intersect the 
section which corresponds to the zero of the universal elliptic curve.
Then the surface is simply a blowup of $\PP^2$ in a single point and 
the action of $SL_2({\bf F}_2)$ comes from the symmetries of a regular 
triangle in the plane. The plane is blown up in the center of the
triangle, and the fibers over the cusps are the proper
preimages of the symmetry lines of this triangle. This description allows
us to calculate the quotient easily. In particular, it has only one
component over the cusp of $SL_2({\bf Z})$ and is smooth in the
neighborhood of the cusp. The cover ${\bf P}^2_1 \rightarrow {\bf
P}^2_1/SL_2({\bf F}_2)$ has 
ramification occurring along the curves which are the fibres over the cusps
of $H/\Gamma(2)$ and the fibres $F_i$ and $F_{\omega}$ 
of the fibration ${\bf P}^2_1 \rightarrow
\bar {H /\Gamma(2)}$ over $SL_2({\bf F}_2)$-orbit of respectively
$i=\sqrt{-1}$ and $\omega (\omega^3=1)$. 
Direct calculation  shows that the action of the generator of 
${\bf Z}_2$ on fibre $F_i$ has two fixed points near which 
it acts as $\bigl({-1 \atop 0} {0 \atop -1} \bigr)$ 
and hence yielding in the quotient two $A_1$-singularities.
On the other hand the action of the generator of ${\bf Z}_3$ on 
$F_{\omega}$ 
has two fixed points of different type. Near one of these points it
acts as $\bigl({\omega \atop 0} {0
  \atop \omega} \bigr)$. Hence the quotient has the 
 singularity which can be resolved by single blow up yielding 
single exceptional curve with self-intersection $-3$.  Near another
point it acts as $\bigl({\omega \atop 0} {0
  \atop \omega^2} \bigr)$, so the singularity is $A_2$.

Moreover, the quotient ${\bf P}^2_1/SL_2({\bf F}_2)$
has a structure of a toric surface, which corresponds
to the fan whose four one-dimensional faces are generated by 
$(1,1)$, $(1,-1)$, $(-1,-1)$ and $(-1,2)$. The structure of the fibration
comes from the projection along $x=y$ line. Notice that the cusp is 
not torus-invariant. The meaning of this toric structure is not clear, 
but it is certainly very useful in calculations of dimensions of 
line bundles over this surface. The easiest way to notice this toric
structure is to employ the the Satake compactification of the quotient of
Siegel modular space of rank 2 by the whole $SP_4(\ZZ)$. The
compactification of $H \times {\bf C}/\Gamma_{1,1}(1)$ is the exceptional
divisor of the blow-up along the infinity
curve. The Satake compactification is known to be the weighted projective
space with weights
$(2,3,5,6)$ (\cite{Igusa}) and the infinity curve corresponds to
$x_5=x_6=0$ for an appropriate choice of $x_5$ and $x_6$. Then the
blow-up is easily calculated by means of toric geometry and the fiber is
the toric surface above. 

One can also obtain the compactification starting form compactification
of $H \times {\bf C}/\Gamma_{1,1}(3)$. 
The latter is ${\bf P}^2_9$ and the
pencil of elliptic curves is a linear pencil containing a cubic curve and
its Hessian (Hesse pencil). Such a pencil has 4 singular fibres
corresponding to four cusps of $\Gamma(3)$ and each singular fibre
consists of three lines formed by sides of a triangle.  Similarly to the
level 2 case $H \times {\bf C}/\Gamma_{1,1}(1)=(H \times {\bf C}/
\Gamma_{1,1}(3))/(\Gamma_{1,1}(1)/\Gamma_{1,1}(3))=H \times {\bf
C}/SL_2({\bf F}_3)$. The quotient of ${\bf P}^2_9$ can be found by
identifying $SL_2({\bf F}_3)$ with the subgroup of $PGL_3({\bf C})$ which
is the image of the quotient of the Shephard-Todd group $G_{25} \subset
GL_3({\bf C})$ yielding the same compactification  as the one obtained 
in using universal level 2 curve.
(**). \footnote{(**) The full Hesse group of order $216$
which is the the latter image of Shephard-Todd group
 acts from this point of view as $\Gamma_{{1 \over 3},{1 \over
  3}}(1)/\Gamma_{1,1}(3)$ which is a semidirect product of $SL_{{\bf
 F}_3}$ and ${\bf F}_3^2$.}

\bigskip {\it 2. Dimensions of the spaces of weak Jacobi forms of zero weight.}
 
Weak Jacobi forms can be identified with sections of line bundles.
In particular, weak Jacobi forms of weight $0$ and index $k$
are sections of the bundle $k[D_{(-1,-1)}]$ where $[D_{(-1,-1)}]$
is the line bundle of the toric divisor which corresponds to the 
point $(-1,-1)$ in the above description of the compactification 
as a toric surface. For small $k$ these dimensions are collected
in the following table.
\begin{equation}
\begin{array}{ccccccccc}
k & \vline   & 0 & 1 & 2 & 3 & 4 & 5 & 6\\ \hline
\dim &\vline & 1 & 1 & 2 & 3 & 4 & 5 & 7
\end{array}
\label{table}
\end{equation}

Alternatively, these dimensions could be derived from the Theorem 9.3
of \cite{EZ}. That theorem states that the ring of weak Jacobi forms
with even weight, any index and no character is $M_*[a,b]$ where
$M_*$ is the ring of modular forms, $a$ has weight $-2$ and index $1$,
while $b$ has weight $0$ and index $1$. Since $M_*$ is generated 
by $e_4$ and $e_6$ of weights $4$ and $6$ respectively, the basis of 
the space of weak Jacobi forms of index $k$ and weight $0$ consists
of $e_4^{\alpha}e_6^{\beta}a^\gamma b^\delta$ where
$\alpha,\beta,\gamma,\delta$ are nonnegative, and
$$
\left\{
\begin{array}{l}
4\alpha+6\beta-2\gamma=0\\
\delta+\gamma=k
\end{array}
\right.
$$
Therefore, the dimension of the space of weak Jacobi forms with index $k$ 
and weight $0$ equals the number of integer solutions of the system of
inequalities
$$
\alpha\geq 0,~\beta\geq0, ~2\alpha+3\beta\leq k.
$$
When $k$ is big, the growth is quadratic, but for small $k\leq 5$ it is linear.
This is responsible for the fact that elliptic genera of Calabi-Yau
varieties 
of small dimension are determined solely by their Hodge numbers.
The following proposition follows easily from the above description.

\begin{prop}
{\rm The ring  of weak Jacobi forms of weight zero is freely generated over
$\CC$ by three elements of degrees $1$, $2$ and $3$ respectively.}
\label{weakring}
\end{prop}

\bigskip {\it 3. Elliptic genera of Calabi-Yau varieties.} 
We will now show that elliptic genera of Calabi-Yau varieties span the
corresponding spaces of weak Jacobi forms. 
\begin{theo}
{\rm
For any $d>0$ elliptic genera of Calabi-Yau varieties of dimension $d$
span the space of  weak Jacobi forms of weight $0$ and index $d/2$ (with
a character if $d$ is odd).
}
\label{span}
\end{theo}

{\em Proof.} Let us first handle the case of even dimension. In view of
Proposition \ref{weakring}, it is enough to show that elliptic genera of
Calabi-Yau varieties of dimensions $2$, $4$ and $6$ span the corresponding
spaces. Moreover, the $q^0$ part of the elliptic genus of $M$ can be
written in terms of Hodge numbers $$Ell(M,y,0)=y^{-\dim M/2}
\sum_{p,q}(-1)^{p+q}h^{p,q}y^p(M)$$ so it is enough to show that the
spaces spanned by above polynomials for all $M$ have expected dimension.
When $\dim M=2$ $K3$ surface provides the necessary example. For dimension
4, we can look at $K3\times K3$ and sextic $X_6$ in $\PP^5$. Finally, in
dimension $6$ we look at $K3\times K3\times K3$, $K3\times X_6$ and
$X_8\subset \PP^7$.

In the case of odd dimension $d$ notice that the space of weak Jacobi
forms of weight $0$ and index $d/2$ and character is isomorphic to the
space of weak Jacobi forms with of weight $0$ and index $(d-3)/2$. Really,
transformation properties imply that all Jacobi forms with character
vanish at $z=1/2$, $z=\tau/2$ and $z=(\tau+1)/2$. On the other hand, there
is a unique (up to a constant) weak Jacobi form with character of index
$3/2$ given by $$f(z,\tau)={{\theta (2z,\tau)} \over {\theta(z.\tau)}}.$$
It appears in \cite{Neumann} in the description of the elliptic genera of
Calabi-Yau threefolds. For any weak Jacobi form $g$ of index $d/2$ the
function $g/f$ is holomorphic on $H\times \CC$. Moreover, it is easy to
see that it has a $q$-expansion, so it is a weak Jacobi form of index
$(d-3)/2$ and no character.

As a result, to show that elliptic genera of Calabi-Yau varieties of odd
dimension span the whole space of Jacobi forms it is enough to consider
products of even-dimensional Calabi-Yau varieties with a quintic in
$\PP^4$ (or any other Calabi-Yau threefold with non-zero Euler
characteristic). \hfill$\Box$

{\it 4. Elliptic genera and Hodge numbers.} The goal of this subsection is
to discuss the relation between elliptic genera and Hodge numbers. As was
already mentioned in the proof of Theorem \ref{span},
$$Ell(M,y,0)=y^{-\dim M/2} \sum_{p,q}(-1)^{p+q}h^{p,q}y^p(M)$$ so to know
$Ell(y,0)$ is to know $\chi^p=\sum (-1)^qh^{p,q}$ for all $p$. On the
other hand, $Ell(y,q)$ certainly can't be used to determine individual
$h^{p,q}$ even in dimension $2$.

\begin{theo} 
{\rm
If dimension of a Calabi-Yau manifold is less than 12 or is
equal to 13, then the numbers $\chi_p$ determine its elliptic genus
uniquely. In all other dimensions there exist Calabi-Yau manifolds with
the same $\{\chi_p\}$ but distinct elliptic genera. Here we allow
disconnected Calabi-Yau manifolds. 
}
\label{nothodge}
\end{theo}

{\em Proof.} For $M$ of even dimension $d$ there are at most $d/2$
independent numbers $\chi_p$. Really, because of the symmetry of Hodge
diagram we only need to consider $\chi_0,...,\chi_{d/2}$ and there is also
an additional linear relation, see \cite{LW}. For $M$ of odd dimension $d$
there are $(d-3)/2$ independent numbers $\chi_p$, because we also have
$\chi_0=0$.

There is a ring homomorphism from the ring spanned by all elliptic genera
to the ring spanned by all $Ell(M,y,0)$.  The table (\ref{table}) shows
that this homomorphism can not be surjective in even dimensions starting
at 12 and in odd dimensions starting at 15. It is easy to see that some
polynomial in the generators of the ring of Jacobi forms (constructed from
$K3$, $X_6$ and $X_8$) with integer coefficients goes to zero when
restricted to $q=0$. Then one collects positive and negative terms and
interprets them as elliptic genera of disjoint unions of products of $K3$,
$X_6$ and $X_8$. This proves the second part of the theorem. To prove the
first part, all we need to show is that the spaces spanned by $\{\chi_p\}$
have dimension $d/2$ for $d=2,4,6,8,10$, which is an explicit calculation.
\hfill{$\Box$}

\begin{rem}
{\rm 
We have first conjectured the above result after extensive Mathematica
calculations of elliptic genera and Hodge numbers in terms of Chern
classes.}
\end{rem}

\section{Elliptic genera of toric varieties}
The goal of this section is to calculate elliptic genus $Ell(y,q)$ of an
arbitrary complete smooth toric variety.  As a byproduct, we notice a
curious identity and give a direct elementary proof of it. One can also
define elliptic genus of an arbitrary complete toric variety with only
Gorenstein singularities but we delay it until Section 7 when
the required techniques are developed.

Let us first describe the data that define a complete toric variety
of dimension $d$, see \cite{Danilov, Fulton, Oda}.  We have a lattice 
(which here simply means a free abelian group) $M$ of rank $d$, its dual
lattice $N$ and a complete polyhedral fan $\Sigma$ in $N$. If all cones
$C^*\in \Sigma$ of maximum dimension are simplicial and are generated by a
basis of $N$, then the corresponding variety $\PP_\Sigma$ is smooth, and
vice versa. We will denote this variety simply by $\PP$.

We have observed in Section 3 that elliptic genus $Ell(y,q)$ of
any smooth variety could be calculated as an Euler
characteristics of a certain double graded sheaf. Here we
change our viewpoint a little bit and work with supersheaves.
For the purposes of this section, supersheaves are simply 
usual sheaves with $\ZZ/2\ZZ$ grading. The corresponding graded 
components are called even and odd.
\begin{prop}
{\rm
Let $X$ be a smooth complete variety of dimension $d$. Consider the
following double graded supersheaf on it which we still denote by
$\Ell$ abusing notations slightly.
$$
\Ell(X)=\otimes_{k\geq 1} \wedge (yq^{k-1} T^*)
\otimes_{k\geq 1} \wedge (y^{-1}q^k T)
\otimes_{k\geq 1} \Sym (q^k T^*)
\otimes_{k\geq 1} \Sym (q^k T), 
$$
where the grading is given by the powers of $y$ and $q$ and 
the parity is given by the parity of the degree of $y$.
Then $Ell(X;y,q)$ is the $y^{-d/2}$ times the super Euler characteristics
of the above sheaf. As usual, super Euler characteristics means the 
Euler characteristics of the even part minus the Euler characteristics
of the odd part.
}
\label{euler}
\end{prop}

{\em Proof.} This result follows directly from Hirzebruch-Riemann-Roch
Theorem and definition of $Ell(y,q)$. $\hfill{\Box}$

When $X=\PP$ is a toric variety, we get some help from the action of 
the torus $(\CC^*)^d$. The idea is to somehow split the calculation
into the sum over the characters $m\in M$. Cohomology of $\Ell(\PP)$
could be calculated as \v{C}ech cohomology for the open affine covering
$\{\AA_C={\rm Spec}(\CC[C])\}$ defined by cones $C^*\in \Sigma$.
Intersection of any number of such subsets is another open subset of this
type, so the covering is acyclic for $\Ell(\PP)$. For fixed powers of
$y$ and $q$ the cohomology of the corresponding graded part of the 
\v{C}ech complex is finite-dimensional, but the entries of the \v{C}ech
complex itself have, of course, infinite dimension. Fortunately,
the action of $(\CC^*)^d$ naturally extends to $\Ell(\PP)$, so sections
of $\Ell(\PP)$ over any such affine subset $\AA_C$ admit natural grading
by the lattice $M$, which is compatible with \v{C}ech differential. The
important remark is that for a fixed $C$, $m$ and powers of $y$ and $q$,
the corresponding graded part of the \v{C}ech complex is
finite-dimensional. Thus, to calculate the Euler characteristics of
$\Ell(\PP)$, one can simply calculate the alternating sum of the
dimensions of these spaces for a fixed $m$ and then add the
results for all $m$. This gives us
$$
Ell(\PP;y,q)=y^{-d/2}\sum_{m\in M}(\sum_{C_0^*,...,C_k^*}
(-1)^{k} \dim_m H^0(\AA_{C_0}\cap ... \cap \AA_{C_k},\Ell(\PP)))
$$
where $\dim$ is understood in a super sense (dimension of even minus
dimension of odd). It is convenient to account for all $m$ simultaneously
by introducing a multi-variable $t$ and using $\sum_{m\in M}t^m\dim_m$. In 
the end we will put $t=1$, but for now we will keep it to assure
convergence of all expressions.

Let us first consider $\AA_C$ with $C$ of maximum dimension. It is
isomorphic to the affine space of dimension $d$ with coordinates
$x^{m_1},...,x^{m_d}$ where $m_1,...,m_d$ generate the cone $C$ and form
the basis of $M$. Sections of $\Ell(\PP)$ over this affine set could be
easily calculated. Really, it is a superpolynomial ring (polynomial ring
of even variables tensored with the exterior algebra of the space of odd
variables) with even variables $x^{m_i}, d(x^{m_i})q^k, k\geq 1,
(\partial/\partial x^{m_i}) q^{k-1}, k\geq 1$ and odd variables
$d(x^{m_i})q^ky,k\geq 1, (\partial/\partial x^{m_i}) q^{k-1}y, k\geq 1$.
Therefore,
$$\sum_{m\in M} t^m \dim_m H^0(\AA_C,\Ell(\PP)) = 
\prod_{i=1,...,d} \prod_{k\geq 1} 
{(1-t^{m_i}yq^{k-1})(1-t^{-m_i}y^{-1}q^k)
\over
(1-t^{m_i}q^{k-1})(1-t^{-m_i}q^k).
}
$$
Unfortunately, in this form it is hard to combine together information
from several cones. We are helped by the following observation.
\begin{prop}
{\rm
$$
\prod_{i=1,...,d} \prod_{k\geq 1}
{(1-t^{m_i}yq^{k-1})(1-t^{-m_i}y^{-1}q^k)
\over
(1-t^{m_i}q^{k-1})(1-t^{-m_i}q^k)
}
=
\sum_{m\in M} t^m
\prod_{i=1,...,d} \left({1\over 1-yq^{m\cdot n_i}}\right)
G(y,q)^d$$
where
$$G(y,q)=\prod_{k\geq 1}
{(1-yq^{k-1})(1-y^{-1}q^k)\over(1-q^k)^2}$$
and $n_i$ are generators of $C^*$.
}
\label{torictrick}
\end{prop}

{\em Proof.} It is sufficient to show this for $d=1$, that is 
we need to show that 
\begin{equation}
\prod_{k\geq 1}
{(1-tyq^{k-1})(1-t^{-1}y^{-1}q^k)
\over
(1-tq^{k-1})(1-t^{-1}q^k)
}
=
\sum_{m\in \ZZ} \left({t^m\over 1-yq^m}\right)
\prod_{k\geq 1}
{(1-yq^{k-1})(1-y^{-1}q^k)\over(1-q^k)^2}.
\label{key}
\end{equation}
We are aware of two proofs of this identity. The first one uses 
some representation theory and is based on the calculation of 
cohomology of a certain explicit operator on the Fock space 
of two free bosons and two free fermions. It is one of the main
ingredients of the preprint \cite{Bvertex}. There also exists
an elementary proof which is sketched below. If we treat
$t$ and $q$ as complex numbers with $|q|<|t|<1$ then both sides
of the equation are meromorphic functions of $z$ (recall that 
$y=\ee^{2\pi\ii z}$). It is straightforward to show that the 
ratio of the two sides is double periodic, because both sides
acquire a factor of $-t^{-1}y^{-1}$ under $y\to yq$. Moreover,
the ratio 
$$
\sum_{m\in \ZZ} \left({t^m\over 1-yq^m}\right)
\prod_{k\geq 1}
{(1-yq^{k-1})(1-y^{-1}q^k)(1-tq^{k-1})(1-t^{-1}q^k)
\over
(1-q^k)^2(1-tyq^{k-1})(1-t^{-1}y^{-1}q^k)
}
$$
has no poles. Really, because of the periodicity, it is enough to
check that there is no pole at $y=t^{-1}$, which
follows from 
$$\sum_{m\in \ZZ}{t^m\over 1-t^{-1}q^m}=
\sum_{m>0}\sum_{k\geq 0} t^m t^{-k}q^{mk}+{1\over 1- t^{-1}}+
\sum_{m>0}\sum_{k\geq 1} (-1) t^{-m} t^{k} q^{mk}
$$
$$
=\sum_{m>0}t^{m} +{1\over1-t^{-1}}+\sum_{n>0}q^n\sum_{mk=n,m,k>0}(t^{m-k}-
t^{k-m})
=0.
$$
Then it remains to observe that at $y=1$ both sides are equal to $1$.
\begin{rem}
{\rm Expressions $1\over 1-yq^{m\cdot n_i}$ in the
above formulas could be understood as power series in $q$
whose coefficients are rational functions in $y$, that is when $m\cdot
n_i<0$, they still have to be expanded around $q=0$. Basically, when
you multiply it by $G(y,q)^d$ as we do, the denominators disappear and 
the resulting functions should be expanded around $q=0$, so we might as
well expand the original expressions around $q=0$.
}
\label{q=0}
\end{rem}

One can also show that the same formula holds for sections of 
$\Ell(\PP)$ over $\AA_C$ where $C^*$ is allowed to have a dimension
smaller than $d$. 
\begin{prop}
{\rm 
For any cone $C^*$ in $\Sigma$ we have 
$$\sum_{m\in M} t^m \dim_m H^0(\AA_C,\Ell(\PP)) =
\sum_{m\in M} t^m
\prod_{i=1,...,\dim C^*} \left({1\over 1-yq^{m\cdot n_i}}\right)
G(y,q)^d$$
where
$$G(y,q)=\prod_{k\geq 1}
{(1-yq^{k-1})(1-y^{-1}q^k)\over(1-q^k)^2}$$
and $n_i$ are generators of $C^*$.
}
\label{anysmoothcone}
\end{prop}

{\em Proof.} The corresponding $\AA_C$ could be (non-canonically) split as
$\CC^{\dim C^*}\times (\CC-\{0\})^{d-\dim C^*}$. Then over the torus we
can use $dx\over x$ and $x\partial_x$ as a basis for differential forms
and vector fields. They have zero grading, and it is then easy to
calculate the contribution of a torus. Combined with the previous
proposition for $\CC^{\dim C^*}$, this gives the desired formula.
$\hfill\Box$

We are now in a position to formulate the main result of this section.
\begin{theo}
{\rm 
Let $\PP$ be a smooth toric variety given by the fan $\Sigma$.
Then
$$
Ell(\PP,y,q)=y^{-d/2}
\sum_{m\in M}
\sum_{C^*\in\Sigma}
(-1)^{{\rm codim} C^*} 
\left(
\prod_{i=1,...,\dim C^*} {1\over 1-yq^{m\cdot n_i}}
\right)
G(y,q)^d$$
where
$$G(y,q)=\prod_{k\geq 1}
{(1-yq^{k-1})(1-y^{-1}q^k)\over(1-q^k)^2}$$
and $n_i$ are generators of $C^*$.
}
\label{elltoric}
\end{theo}

{\em Proof.}
Because of Proposition \ref{anysmoothcone} it is enough to show that
in the alternating sum over the \v{C}ech complex each cone
$C^*\in\Sigma$ is counted with the coefficient $(-1)^{{\rm codim} C^*}$.
We will carry an induction on codimension of $C^*$.
Obviously, the cones that do not contain $C^*$ could be ignored. 
Then we can quotient out the subspace generated by $C^*$. The cones that
contain $C^*$ will form a complete fan on the remaining subspace of 
dimension ${\rm codim} C^*$. The sum of $(-1)^k$ over all
possible intersections of $k+1$ cones in this fan is, of course, $1$
by binomial formula. On the other hand this sum splits according to 
possible intersections of the cones in the fan. All cones, except
the vertex (image of $C^*$) contribute $(-1)^{\rm codim}$ by the induction
assumption. So the occurrence of $C^*$ is 
$$1-\sum_{C_1*,C^*\subset C_1^*} (-1)^{{\rm codim} C_1^*}$$
and the above sum could be expressed as the Euler characteristics of 
a sphere of dimension ${\rm codim}C^*-1$, which finishes the proof.
$\hfill{\Box}$

\begin{rem} 
{\rm Alternatively, one could use a slight modification of the
\v{C}ech complex where you use each cone once and the differential
consists of the restriction maps with the signs as in the singular
homology complex of $\Sigma$. Then the above formula is automatic, 
but the difficult part is to show that the complex can be used to
calculate the cohomology of any coherent sheaf. 
}
\end{rem}

Since a complete toric variety can never be a Calabi-Yau, we are mostly
interested here in the $y=-1$ case. Also, one can take $G(-1,q)^d$ outside
of the summation sum once we agree to expand around $q=0$, see Remark
\ref{q=0}. If we denote 
$$\widehat{Ell}(X;q)=(-1)^{d/2}
Ell(X;-1,q) G(-1,q)^{-d}$$ we get the following result. 
\begin{theo}
{\rm If $\PP$ is a smooth complete toric variety then 
$$\widehat{Ell}(\PP;q)=\sum_{m\in M}
\left(
\sum_{C^*\in\Sigma}
(-1)^{{\rm codim} C^*}
\prod_{i=1,...,\dim C^*} {1\over 1+q^{m\cdot n_i}}
\right).
$$
The inner sum here is taken first. Then for a given degree of $q$
only a finite number of $m$ contribute to $\widehat{Ell}$.
}
\end{theo}

{\em Proof.} The only thing to prove here is the last statement.
However, the part of the sheaf $\Ell(\PP)$ for a fixed degree of $q$ 
is coherent, so its cohomology has finite dimension. Therefore, only 
a finite number of $m$ contribute to $Ell(\PP,-1,q)$ at a given degree.
Multiplication by $G(-1,q)^{-d}$ does not change this.
$\hfill{\Box}$

\begin{rem}
{\rm
Notice that $\widehat{Ell}(\PP)$ is precisely the Landweber-Stong genus,
see \cite{Landweber.Stong}. In particular, it is known to be a
modular form with respect to the index three subgroup $\Gamma_0(2)$ in
the full modular group. However, we do not see it from the formulas above.
One may also wonder if these formulas could be modified to give examples
of modular forms at different levels. An obvious thing to try is to
use $\sum_{m\in M'}$ where $M'$ is a sublattice of $M$.
}
\end{rem}

\begin{rem} 
{\rm 
We can compare our result for $\PP=\PP^2$ with the known
Landweber-Stong genus of $\PP^2$. After rather easy simplifications,
this results in the following curious identity.
$$
\sum_{m\geq 1,n\geq 1} {q^{m+n}\over (1+q^m)(1+q^n)(1+q^{m+n})}=
\sum_{r\geq 1}q^{2r}\sum_{k|r}k.
$$
}
\end{rem}

Even though our discussion above provides another proof of it, there
exists a direct derivation. Also, George Andrews have shown us how
to reduce this identity to the result of \cite{Bailey}.

{\em Proof of the $\PP^2$ identity.}
It is easy to show that
$${xy\over (1+x)(1+y)(1+xy)}=
\sum_{a\geq1, b\geq 1, {\rm min}(a,b)=odd}x^ay^b(-1)^{a+b},$$
which implies that
$$
\sum_{m\geq 1,n\geq 1} {q^{m+n}\over (1+q^m)(1+q^n)(1+q^{m+n})}=
\sum_{d\geq 2}q^d
\sum_{
a\geq 1, b\geq 1, m\geq 1, n\geq 1,
ma+nb=d,{\rm min}(a,b)=odd}
\hspace{-8pt}(-1)^{a+b}.
$$
For a given $d$ we look at the set of all solutions to
$a\geq 1, b\geq 1, m\geq 1, n\geq 1,
ma+nb=d,{\rm min}(a,b)=odd$ which we denote by $I$. We denote by $J$
the set obtained from $I$ by excluding all solutions with $m=n$
and $a+b=even$. We will show that the set $J$ does not contribute to
the above sum. We denote by $J_{odd}$ and $J_{even}$ the parts of $J$
with odd and even $a+b$ respectively and plan to show that there is
a one-to-one correspondence between these two sets.
  
Given an element $(a,b,m,n)$ in $J_{even}$ we construct an element in
$J_{odd}$ either as $(a+b,b,m,n-m)$ or $(a,a+b,m-n,n)$ depending on
which of $m$ and $n$ is bigger (they can not be equal, because we are in
$J$). Given an element $(a,b,m,n)$ in $J_{odd}$ we construct an
element in $J_{even}$ as $(a-b,b,m,m+n)$ or $(a,b-a,m+n,n)$ depending
on which of $a$ and $b$ is bigger. One can show that these two maps are
well-defined and are inverses of each other.

So we have
$$
\sum_{m\geq 1,n\geq 1} {q^{m+n}\over (1+q^m)(1+q^n)(1+q^{m+n})}
=
\sum_{d\geq 2} q^d \sum_{a\geq 1,b\geq 1, m\geq 1,m(a+b)=d,
{\rm min}(a,b)=odd, a+b=even} 1
$$
which is easily seen to equal 
$$
\sum_{r\geq 1}q^{2r}\sum_{k|r}k. 
$$\\[-3em]
$~$$\hfill{\Box}$

\section{Elliptic genera and mirror symmetry} 
The goal of this section is to show that elliptic genera of two mirror
symmetric hypersurfaces in toric varieties coincide up to an expected 
sign. Because of Theorem \ref{nothodge}, this provides a new check
of mirror symmetry in high dimension. We introduce a concept of
elliptic genus for hypersurfaces in toric Fano varieties with Gorenstein 
singularities, and prove the mirror duality for these genera as well.
Unfortunately, the argument uses full force of the machinery of
\cite{Bvertex} and \cite{MSV}, even in the smooth case. We will be reviewing
briefly necessary techniques as we go on.

Chiral de Rham complex $\msv(X)$ of any smooth variety $X$ was defined in
\cite{MSV}. It is a sheaf of vector spaces over $X$ equipped with 
a double grading by eigenvalues of $J[0]$ and $L[0]$. It is also given a
$\ZZ/2\ZZ$ (even-odd) grading induced from the $Z$-grading by $J[0]$. 
We must mention that ${\cal MSV}(X)$ itself in not a 
quasi-coherent sheaf. The "multiplication by functions" map of sheaves
${\cal O}(X)\times{\cal MSV}(X)\to{\cal MSV(X)}$ is defined but is 
not associative. However, ${\cal MSV}(X)$ possesses a natural filtration
compatible with the grading and the above multiplication such that the
graded object is a quasi-coherent sheaf isomorphic to
$$
\Ell=\otimes_{k\geq 1} \wedge (yq^{k-1} T^*)
\otimes_{k\geq 1} \wedge (y^{-1}q^k T)
\otimes_{k\geq 1} \Sym (q^k T^*)
\otimes_{k\geq 1} \Sym (q^k T),
$$
see \cite{MSV}. Here the powers of $y$ and $q$ stand for the eigenvalues
of $J[0]$ and $L[0]$ respectively. Since at each power of $q$ in the above
expression we have a coherent sheaf (this means that $\cal MSV$ is
loop-coherent, or loco, in the terminology of \cite{Bvertex}), the
cohomology of ${\cal MSV}(X)$ are finite-dimensional vector spaces over
$\CC$ at every given eigenvalue of $L[0]$. As a result, the (double
graded, super) Euler characteristics of $\msv(X)$ is well-defined and
coincides with the Euler characteristics of the sheaf above.  Together
with Proposition \ref{euler}, this prompts the following definition.
 
\begin{dfn}
{\rm
Let $X$ be any variety of dimension $d$ for which there is defined the
\MSV\ $\msv(X)$. Then elliptic genus $Ell(X;y,q)$ is defined by
$$
Ell(X;y,q) = y^{-d/2} {\rm SuperTrace}_{H^*(\msv(X))} y^{J[0]}
q^{L[0]}.
$$
}
\label{EGasTrace}
\end{dfn}
This allows us to extend the notion of elliptic genus to some singular
varieties for which chiral de Rham complexes have been constructed. At
this point, such varieties include arbitrary Gorenstein toric varieties
and Calabi-Yau hypersurfaces in Gorenstein toric Fano varieties, see
\cite{Bvertex}. The above formula makes sense because it was shown in
\cite{Bvertex} that for fixed $(J[0],L[0])$ eigenvalues the corresponding
eigenspace in $H^*$ is finite-dimensional. We have shown in Section 3
that for any smooth Calabi-Yau variety thus defined elliptic genus is a
weak Jacobi form.

We will now describe how to calculate elliptic genera of Calabi-Yau
hypersurfaces in Gorenstein toric Fano varieties in terms of combinatorial
structures that define them. We recall that such a family of Calabi-Yau
hypersurfaces is determined by the following data, see \cite{bat.dual}.

Let $M_1$ and $N_1$ be dual free abelian groups of rank $d+1$. Denote
by $M$ and $N$ two dual free abelian groups such that $M=M_1\oplus {\bf
Z}$ and $N=N_1\oplus{\bf Z}$. Element $(0,1)\in M$ is denoted by $\rm deg$
and element $(0,1)\in N$ is denoted by $\deg^*$. There are dual reflexive
polytopes $\Delta\in M_1$ and $\Delta^*\in N_1$ which give rise to dual
cones $K\subset M$ and $K^*\subset N$. Namely, $K$ is a cone over
$(\Delta,1)$ with vertex at $(0,0)_M$, and similarly for $K^*$. 
There is a complete fan on $N_1$ whose one-dimensional cones are generated
by some lattice points in $\Delta^*$ (in particular, by all vertices).
This fan induces the decomposition of the cone $K^*$ into subcones, each
of which includes $\deg^*$. Let us denote this decomposition by $\Sigma$.
It will later turn out that elliptic genus does not depend on the choice
of $\Sigma$ (that is, partial crepant toric desingularizations do not
alter it).

The formula for the elliptic genus that is the cornerstone of the argument
is given by the following proposition.
\begin{prop}
{\rm
Let $X$ be a generic hypersurface in the Gorenstein toric Fano variety
defined by the combinatorial data above. Then
$$
Ell(X;y,q)=y^{-d/2}\sum_{m\in M} 
\left(
\sum_{n\in K^*}
y^{n\cdot {\deg} -m\cdot {\deg}^*}
q^{m\cdot n+m\cdot{\deg}^*}
G(y,q)^{d+2}
\right)
$$
where 
$$G(y,q)=\prod_{k\geq 1}
{(1-yq^{k-1})(1-y^{-1}q^k)\over(1-q^k)^2}.$$
The inner summation times $G^{d+2}$ is a well-defined double series in $y$
and $q$, when $G$ is expanded around $q=0$. Then the outer sum turns out
to make sense as a double power series in $y$ and $q$, that is only a
finite number of $m$ contribute to a given coefficient. 
} 
\label{formula}
\end{prop}

\begin{rem}
{\rm
This formula clearly shows that elliptic genus does not
depend on the choice of desingularization. In particular, if $X$ admits
a crepant toric desingularization $\hat X$, then
$$Ell(X;y,q)=Ell(\hat X;y,q).$$
}
\end{rem}

Before we begin to prove Proposition \ref{formula}, we need to recall
more background material from \cite{Bvertex}. That
paper contains the calculation of the cohomology of the \MSV\ for toric
varieties and hypersurfaces in them. If the ambient toric variety is
smooth, then $H^*(\msv)$ of the hypersurface can be calculated as the
cohomology of a certain vector space by a certain operator, both described
below.

This space is built from tensor products of irreducible modules of
infinite Heisenberg and Clifford algebra. Those are the affinization 
(resp. Clifford affinization, cf. \cite{Kac} p.26) of 
abelian Lie algebra $M \oplus N$ in which all elements are considered 
even (resp. odd) with supersymmetric (resp. skew supersymmetric) 
bilinear form given by natural pairing $M \times N \rightarrow {\bf Z}$. 
In other words these algebras (resp. $H(M \oplus N)$ and $Cl(M \oplus N)$) 
are 
 $$ (M \oplus N) \otimes {\bf C}[t,t^{-1}]\oplus {\bf C}K$$ 
with commutator $$\{x \otimes t^k, y \otimes
t^l\}_{-}=k(x,y)\delta_{k+l,0}K$$
resp. 
$$ \{x \otimes t^k, y \otimes t^l\}_{+}=(x,y)\delta _{k+l,0} K$$
(here  $\{. ,.\}_-$ and   $\{. ,.\}_+$) denote  
commutator and anticommutator respectively). 
We will always assume that the commutator $K$ acts as identity 
and we identify $K$ for both algebras. We will thus consider
$${\bf g}= (M \oplus N)_H \otimes {\bf C}[t,t^{-1}]\oplus 
(M \oplus N)_{Cl} \otimes {\bf C}[t,t^{-1}]\oplus 
{\bf C}K$$
where subscripts are used to distinguish between two copies of the
lattice.

The ``canonical representations'' of ${\bf g}$ (the Fock spaces) denoted
${\rm Fock}_{m_0 \oplus n_0}$ can be obtained as follows. We define $b{\bf
g}$ as 
$$b{\bf g}=bH(M\oplus N)\oplus bCl(M\oplus N)\oplus\CC K$$
$$ = 
\oplus_{k \geq  0} (M \oplus N)_H \otimes t^k
 \oplus(\oplus_{k > 0} M_{Cl} \otimes t^k)\oplus (\oplus_{k\geq 0}
N_{Cl} \otimes t^k) \oplus {\bf C} K$$ 
and view $\bf C$ as a $b{\bf g}$-module such that 
$$bCl(M\oplus N)\cdot 1=0, K \cdot 1 =1, x_H \otimes t^0 \cdot 1=x
\cdot (m_0 \oplus n_0),$$
$$
 x_H \otimes t^k \cdot 1 =0 \ {\rm for} \ k \geq 1 \ (x \in M \oplus N). 
$$
Then we put 
$${\rm Fock}_{m_0 \oplus n_0}=U({\bf g})\otimes_{b{\bf g}}\CC.$$

We will be using the notations of \cite{Bvertex}. There, algebra
${\bf g}$ is described in terms of generators
$n\cdot A[k], m\cdot B[k], m\cdot \Phi[k], n\cdot\Psi[k]$ for all $k\in \ZZ$
where $n\cdot A[k]$ stands for $(0\oplus n)_H\otimes t^k$ and similarly
for $B,
\Phi,\Psi$.
The only non-trivial super-commutators are
$$\{m\cdot B[k], n\cdot A[l]\}_-=(m\cdot n)k\delta_{k+l}^0,~
\{m\cdot \Phi[k], n\cdot \Psi[l]\}_+=(m\cdot n)\delta_{k+l}^0.$$

As a a vector space ${\rm
Fock}_{m_0\oplus n_0}$ is a tensor product of the symmetric algebra of a
certain infinite-dimensional space and an exterior algebra of another
infinite-dimensional space. The first space is
$$ \oplus_{k\geq 1} M\cdot B[-k]{\bf C} \oplus_{k\geq 1}
N\cdot A[-k]{\bf C}
$$
and the second space is 
$$
\oplus_{k\geq 1} M\cdot \Phi[-k]{\bf C} \oplus_{k\geq 0}
N\cdot \Psi[-k]{\bf C}.
$$
Here $M\cdot B[-k]{\bf C}$ means simply a copy of $M_\CC$, one for each
$k$, and similarly for $A$, $\Phi$ and $\Psi$.

Elements of ${\rm Fock}_{m\oplus n}$ can be written as 
$$ \sum 
(m_{11}\cdot B[-1])(m_{12}\cdot B[-1])...
(m_{21}\cdot B[-2])(m_{22}\cdot B[-2])...
$$
$$
(n_{11}\cdot A[-1])(n_{12}\cdot A[-1])...
(n_{21}\cdot A[-2])(n_{22}\cdot A[-2])...
$$
$$
(m_{11}'\cdot \Phi[-1])(m_{12}'\cdot \Phi[-1])...
(m_{21}'\cdot \Phi[-2])(m_{22}'\cdot \Phi[-2])...
$$
\begin{equation}
(n_{11}'\cdot \Psi[0])(n_{12}'\cdot \Psi[0])...
(n_{21}'\cdot \Psi[-1])(n_{22}'\cdot \Psi[-1])...
|m,n>
\label{monom}
\end{equation}
where there are only finitely many summands and finitely many factors in
each of the summands. It should be understood that the above expressions
are linear in $m_{ij}$ and $n_{ij}$, and $A$, $B$, $\Phi$ and $\Psi$
super-commute. Commutation rule is necessary to understand the relation to
the product of symmetric and exterior powers. 

Operators $m\cdot B[k], k<0$, $n\cdot A[k], k<0$, $m\cdot \Phi[k],k\leq 0$
and $n\cdot \Psi[k],k< 0$ act on this space via multiplication in the 
(super)-polynomial rings. Operators $n\cdot A[0]$ and $m\cdot B[0]$ act
as scalars $m_0\cdot n$ and $m\cdot n_0$ respectively. Other operators
act as derivations on these rings. Our notations are set up in a way that $A[k]$,
$B[k]$, $\Phi[k]$ and $\Psi[k]$ are vector-valued operators on this Fock space.

For any subset of $M\oplus N$ we can consider the direct sum of the
corresponding Fock spaces. For example, we will consider 
$${\rm Fock}_{M\oplus N}=\oplus_{m_0\in M,n_0\in N}{\rm Fock}_{m_0\oplus n_0}$$
and 
$${\rm Fock}_{M\oplus K^*}=\oplus_{m_0\in M,n_0\in K^*}{\rm Fock}_{m_0\oplus
n_0}.$$

Out of the operators $m\cdot B[l]$, $n\cdot A[l]$,
$m\cdot \Phi[l]$ and $n\cdot \Psi[l]$  one can construct the operators $L[0]$ and
$J[0]$ which act on each of the spaces $\Fock_{m_0\oplus n_0}$. The
details could be found in \cite{Bvertex}, but it is useful to mention that the
"monomials" in the polynomial ring written above are eigen-vectors of these
operators.  The eigen-values are calculated as follows.

{\it
$L[0]$ counts the opposite of the sum of indices [ ] in the above expression
(\ref{monom})
of the element, plus $m\cdot n+{\deg}^*\cdot m$;}

{\it $J[0]$ counts the number of occurrences of $\Phi$ minus the number
of 
occurrences of $\Psi$ plus ${\deg}\cdot n - {\deg}^*\cdot m$.}

In addition to these operators there are defined {\em vertex operators}\
${\rm e}^{\int m\cdot A+n\cdot B}[l]$ for all $m$ and $n$ in the
corresponding lattices and all integer $l$. These operators map 
between Fock spaces with different eigenvalues of $A[0]$ and $B[0]$. We
refer to \cite{Bvertex} for the explicit formula. When there is a fan
$\Sigma$ present, one can use it to modify the definition of the
vertex operators to make their results zero unless both $n$ and
$B[0]$-eigenvalue of the argument lie in the same cone. This introduces a
structure of the vertex algebra (see \cite{Kac}) on ${\rm Fock}_{M\oplus
N}$ and some subspaces of it. Even though the underlying vector spaces
stays the same, they are denoted by ${\rm Fock}_{...}^{\Sigma}$ to show
that vertex operators act differently.

A choice of a particular hypersurface $X$ in a family means a choice of 
coefficients $f_m$ for each $m\in (\Delta,1)$. We will be concerned with
a generic such choice. Also, we choose generic set of numbers $g_n$,
one for each $n\in (\Delta^*,1)$. When the ambient variety is smooth,
$g$ plays only a token role. In the presence of singularities the role of
$g$ is more essential, they code for some strange structure that helps
define the chiral de Rham complex. In \cite{Bvertex} there are introduced  
operators ${\rm BRST}_f$ and ${\rm BRST}_g$ 
$$
{\rm BRST}_f=\left(\sum_{m\in (\Delta,1)} (m\cdot \Phi){\rm e}^{\int
m\cdot B}\right)[0]
$$
$$
{\rm BRST}_g=\left(\sum_{n\in (\Delta^*,1)} (n\cdot \Psi){\rm e}^{\int
n\cdot A}\right)[0]
$$
which act on 
${\rm Fock}_{M\oplus N}$, ${\rm Fock}_{M\oplus N}^\Sigma$, 
${\rm Fock}_{M\oplus K^*}$ and ${\rm Fock}_{M\oplus K^*}^\Sigma$.
Moreover, for any cone $C^*\in \Sigma$ these operators act on
${\rm Fock}_{M\oplus C^*}^\Sigma$. Both operators act as differentials
of a double complex, that is they anticommute, their squares are zero,
and they change eigen-values of $({\deg}^*\cdot A[0],{\deg}\cdot
B[0])$ by $(1,0)$ and $(0,1)$ respectively. Besides, they both commute 
with $J[0]$ and $L[0]$ defined above, so the double grading on $\Fock$
descends to the cohomology of these operators.

The following statement defines \MSV\ of a hypersurface in a 
Gorenstein  toric Fano variety.

\begin{dfn}
{\rm
(\cite{Bvertex})
Sections of \MSV\ over the
Zariski open affine chart of the hypersurface that corresponds to the
cone $C^*$ are defined as elements of cohomology of ${\rm Fock}_{M\oplus
C^*}$ with respect to the operator 
$$
{\rm BRST}_{f,g}={\rm BRST}_f+{\rm BRST}_g.
$$
}
\end{dfn}

When the ambient toric variety is nonsingular this definition coincides
with the definition of Malikov, Schechtman and Vaintrob for an arbitrary 
choice of nonzero $g_n$. This construction behaves well under
the localization, and it allows us to define $\msv(X)$ as a quasi-loco
sheaf of vertex algebras, see \cite{Bvertex}.

One of the main results of \cite{Bvertex} is the following. 
\begin{theo}
{\rm (\cite{Bvertex})
If the ambient variety is smooth, cohomology $H^*(\msv(X))$ is 
isomorphic to the cohomology of $\Fock_{M\oplus K^*}^\Sigma$ with respect 
to the operator ${\rm BRST}_{f,g}$. 
}
\label{smoothmsv}
\end{theo}
Unfortunately, it is unclear whether this theorem holds in the general
case. However, we will show below that the (graded, super) Euler
characteristics of $H^*(\msv(X))$ equals the Euler characteristics of 
the above space.

{\em Proof of the Proposition \ref{formula}.}
General theory of quasi-loco sheaves developed in \cite{Bvertex} implies that
the cohomology of $\msv(X)$ can be calculated by means of the  \v{C}ech complex
defined by the toric affine charts. Each entry of the \v{C}ech complex is  
the direct sum of global sections over affine charts that correspond to
cones $C^*$ in $\Sigma$ of various dimensions. This naturally leads to
considering the following double
complex. It has entries labeled by eigenvalues of 
${\deg}^*\cdot A[0]+
{\deg}\cdot B[0]$
 and the position in the \v{C}ech complex. The
differentials are ${\rm BRST}_{f,g}$ and $d_{\hat{C}ech}$. After
some insertions of $(-1)$, these differentials anticommute. 

Cohomology of $\msv(X)$ is defined as the repeated cohomology of this 
complex when you first take cohomology with respect to 
${\rm BRST}_{f,g}$ and then use \v{C}ech differential.
On the other hand, one can start with cohomology of \v{C}ech differential
and then do ${\rm BRST}_{f,g}$. It is easy to show that cohomology of
the  \v{C}ech differential is non-zero only at the zeroth column, where 
it is isomorphic to ${\rm Fock}_{M\oplus K^*}^\Sigma$, and the cohomology 
of the total complex is the cohomology of ${\rm Fock}_{M\oplus K^*}$ with
respect to ${\rm BRST}_{f,g}$. 
Therefore, there exists a spectral sequence from $H^*(\msv(X))$ to this
space. If the ambient variety is smooth, then it is proved in \cite{Bvertex} to
degenerate immediately, providing an explicit description of the cohomology of
$\msv(X)$, see Theorem \ref{smoothmsv}. Since the differentials of this
spectral sequence change parity and commute with $L[0]$ and $J[0]$, they
will have no effect on the supertrace used to define $Ell(y,q)$. So we
have
$$
Ell(y,q)=y^{-d/2}{\rm SuperTrace}_{H^*{\rm Fock}_{M\oplus K^*}^\Sigma}
y^{J[0]}q^{L[0]}
$$
where $H^*$ denotes cohomology with respect to ${\rm BRST}_{f,g}$. 
We have used the fact that each of the bigraded components of $\msv(X)$ 
has finite-dimensional cohomology.

To calculate the supertrace above, we consider another double
complex where the entries are $({\deg}^*\cdot A[0],{\deg}\cdot
B[0])$ eigenspaces of ${\rm Fock}_{M\oplus K^*}^\Sigma$ and the 
differentials are ${\rm BRST}_f$ and ${\rm BRST}_g$. Eigenvalues
of ${\deg}\cdot B[0]$ are non-negative, because ${\deg}\cdot K^*\geq
0$, but eigenvalues of ${\deg}^*\cdot A[0]$ are not. The supertrace
over the cohomology of the double complex equals to the supertrace over
the repeated cohomology if we are able to show that the corresponding
spectral sequence converges. Notice that we can split the double complex
according to eigenvalues of $L[0]$ and $J[0]$ because these operators
commute with differentials. We take cohomology with respect to ${\rm
BRST}_g$ first. One can show that this cohomology is zero for all
sufficiently big values of ${\deg}\cdot B[0]$, see 
\cite{Bvertex}, which implies that the spectral sequence degenerates
after a finite number of steps. Moreover, for a fixed pair of
eigen-values of $L[0]$ and $J[0]$ all vertical cohomology spaces are
finite-dimensional (we will show it below), so we can
calculate the super-trace on the first layer of the spectral sequence.

\begin{lem}
{\rm
For a fixed pair of eigenvalues of $L[0]$ and $J[0]$ the
corresponding eigen-space of ${\rm BRST}_g$ cohomology of 
${\rm Fock}_{M\oplus K^*}^\Sigma $ or 
${\rm Fock}_{M\oplus K^*}$ 
is finite-dimensional. 
}
\label{findim}
\end{lem}

{\em Proof of the lemma.} Consider the grading operator 
$L[0]+J[0]$. It counts the opposite of the sum of mode numbers
plus the number of $\Phi$ minus the number of $\Psi$ plus
${\deg}\cdot n$. 

If we denote by $\PP$ and $L$ the ambient toric
variety and the canonical line bundle over it, then cohomology
of $\msv(L)$ is precisely the ${\rm BRST}_g$ cohomology of
${\rm Fock}_{M\oplus K^*}^\Sigma $. It was proved in 
\cite{Bvertex} that $\msv(L)$ is a loco sheaf 
with respect to the above grading operator $L[0]$. Therefore,
its cohomology has a filtration such that all the quotients 
are finitely generated modules over $\CC[K]$, where the action
of $\CC[K]$ is given by 
$\ee^{\int m\cdot B}[0]$. However,
a multiplication by $\ee^{\int m\cdot B}[0]$ decreases the eigenvalue
of $J[0]$ by ${\deg}^*\cdot m$, which shows that for a fixed
value of $L[0]$ and $J[0]$ the eigen-space is
finite-dimensional. A case of ${\rm Fock}_{M\oplus K^*}$ is treated 
analogously, because it admits a filtration whose quotients are finitely 
generated  $\CC[K]$ modules.
$\hfill{\Box}$

So now we can take supertrace over the cohomology of ${\rm Fock}_{M\oplus
K^*}^\Sigma$ by ${\rm BRST}_g$. We notice that this differential commutes
with $A[0]$ so we can split the whole picture according to eigenvalues of 
$A[0]$. Basically, we have managed to move the problem from the Calabi-Yau
hypersurface to the line bundle over the ambient toric variety which 
has a huge torus symmetry. 

Let $m\in M$ be one such eigenvalue. We claim that for fixed values of
$m$, $J[0]$ and $L[0]$ the dimension of the corresponding eigen-space of
${\rm Fock}_{M\oplus K^*}^\Sigma$ is finite. Really, look at the
$L[0]+rJ[0]$ eigenvalue for sufficiently big fixed $r$, which is chosen in
a way that $m\cdot n+(r-1){\deg}\cdot n\geq 0$ for all $n\in K^*$. For
any element $|m,n>$ this eigenvalue is bounded from below, and only
finitely many $n$ work for each particular value. Moreover, multiplying by
$A,B,\Phi,\Psi[-l]$ only increases this value, except for a finite number
of anticommuting modes $\Psi[-l]$, which can decrease the eigenvalue by no
more than a constant. 

As a result, for a fixed $m$ we can calculate the supertrace over 
the ${\rm BRST}_g$ cohomology of ${\rm Fock}_{m\oplus K^*}$ directly
over the Fock space. We get 
$${\rm SuperTrace}_{{\rm Fock}_{m\oplus K^*}}y^{J[0]}q^{L[0]}=
{\rm SuperTrace}_{{\rm Fock}_{m\oplus 
K^*}}
(yq^{-r})^{J[0]}q^{L[0]+rJ[0]}=$$
$$
= 
\left(\sum_{n\in K^*}
y^{n\cdot {\deg} -m\cdot {\deg}^*}
q^{m\cdot n+m\cdot{\deg}^*}\right)
G(y,q)^{d+2}
$$
where
$$G(y,q)=\prod_{k\geq 1}
{(1-yq^{k-1})(1-y^{-1}q^k)\over(1-q^k)^2}.$$
Really, for every $n\in K^*$ the space $\Fock_{m\oplus n}$ is 
isomorphic to the tensor product of spaces 
$$\CC\oplus m_i\cdot \Phi[-k]\CC, k\geq 0$$
$$\CC\oplus n_i\cdot \Psi[-k]\CC, k>0$$
$$\CC\oplus m_i\cdot B[-k]\CC\oplus (m_i\cdot B[-k])^2\CC\oplus ...,
~k> 0$$
$$\CC\oplus n_i\cdot A[-k]\CC\oplus (n_i\cdot A[-k])^2\CC\oplus ...,
~k>0$$
where $i=1..d$ and $\{m_i\}$ and $\{n_i\}$ are arbitrary bases of 
$M_\CC$ and $N_\CC$. When one calculates the supertrace, the first 
two spaces contribute to the numerator and the last two to the denominator
of $G(y,q)$. The term 
$$y^{n\cdot {\deg} -m\cdot {\deg}^*}
q^{m\cdot n+m\cdot{\deg}^*}$$
appears when one evaluates $J[0]$ and $L[0]$ on $|m,n\vac$ itself.

We remark that it is legal to multiply these two double series above
because the former has only a finite number of entries for each
$L[0]+rJ[0]$ eigenvalue and coefficients of the second one are supported
inside a parabola. Moreover, our arguments show that the resulting double
series in $y$ and $q$ will contain entries of a fixed bi-degree only for a
finite number of $m$.

When we add the supertraces above for all $m$ we recover the formula 
of Proposition \ref{formula}, which finishes its proof. $\hfill{\Box}$

Our goal now is to compare the elliptic genera for dual Calabi-Yau
hypersurfaces. This means switching $M, K$ and $N, K^*$.
One would expect to find a direct argument based simply on the 
formula of Proposition \ref{formula}, but we were unable to find one.
The key result of \cite{Bvertex} which we will use here is
\begin{theo}
{\rm
(\cite{Bvertex})
Cohomology of $\Fock_{M\oplus K^*}$ with respect to ${\rm BRST}_{f,g}$ 
is isomorphic to the cohomology of $\Fock_{K\oplus N}$ with respect
to ${\rm BRST}_{f,g}$.

}
\label{keyswitch}
\end{theo}
We remark that this statement is by no means obvious. Also, 
$\Fock_{K\oplus N}$ is not exactly the analog of $\Fock_{M\oplus K^*}$
for the dual pair. The difference comes from the fact that 
in the formula (\ref{monom}) only negative modes of $\Psi$ are allowed,
as opposed to all non-positive modes of $\Phi$. In the mirror 
picture non-positive modes of $\Psi$ and negative modes of $\Phi$
appear. However, one can still construct an isomorphism between
$\Fock_{K\oplus N}$ and the mirror analog of $\Fock_{M\oplus K^*}$ by
shifting the mode numbers of $\Phi$ and $\Psi$ by $-1$ and $1$ 
respectively. Under this isomorphism, $L[0]$ and $J[0]$ are related
to their mirror versions as follows
\begin{equation}
L_{X}[0]= L_{X^*}[0]+J_{X^*}[0],~
J_{X}[0]=-J_{X^*}[0].
\label{newlq}
\end{equation}
This is easy to check from the description of the action of $J[0]$
and $L[0]$ on the monomials in (\ref{monom}).

\begin{prop}
{\rm 
Let $X$ and $X^*$ be two mirror hypersurfaces in Gorenstein toric
Fano varieties. Then 
$$Ell(X;y,q)=y^{-d} q^{d/2} Ell(X^*;y^{-1}q,q).$$
}
\label{duality}
\end{prop}

{\em Proof.} 
The key observation here is that the supertraces calculated over
the BRST cohomology of ${\rm Fock}_{M\oplus K^*}^\Sigma$ and 
${\rm Fock}_{M\oplus K^*}$ are the same. Really, we can prove the same
supertrace formula for the ${\rm BRST}_{f,g}$ cohomology of 
${\rm Fock}_{M\oplus K^*}$. On the other hand, Theorem \ref{keyswitch}
shows that the cohomology of ${\rm Fock}_{M\oplus K^*}$
is isomorphic as a vector space to the cohomology of ${\rm Fock}_{K\oplus
N}$ for the same map. It was important here that we got rid of the 
fans. Now we use formula (\ref{newlq}) to get
$$Ell(X;y,q)=y^{-d/2}{\rm SuperTrace}_{H^*\Fock_{M\oplus K^*}}y^{J_X[0]}
q^{L_X[0]}
$$
$$
=y^{-d/2}{\rm SuperTrace}_{H^*\Fock_{K\oplus N}} y^{-J_{X^*}[0]}
q^{L_{X^*}[0]+J_{X^*}[0]}
= y^{-d} q^{d/2} Ell(X^*;y^{-1}q,q). 
$$\\[-2em]
$~\hfill{\Box}$

\begin{rem}
{\rm
When $X$ (or $X^*$) is smooth, we can use the fact that $Ell(y,q)$
is a weak Jacobi form to show that
$$Ell(X;y,q)=(-1)^d Ell(X^*;y,q)$$
as predicted by Mirror Symmetry. 
}
\end{rem}

\section{Modular properties of elliptic genera of singular varieties}

We will now show that $Ell(X;y,q)$ is a weak Jacobi form for an 
arbitrary reflexive polytope $\Delta$. This will allow us to prove mirror
duality of elliptic genera of hypersurfaces in full generality.
In addition we will extend the definition of elliptic genera for toric
varieties to the Gorenstein case and prove their transformation
properties under $\Gamma_0(2)$.

\begin{theo}
{\rm
Elliptic genus of a generic Calabi-Yau hypersurface of dimension $d$
in any toric Gorenstein Fano variety is a weak Jacobi form of weight 0
and index $d/2$ (with a character if $d$ is odd).
}
\label{modular}
\end{theo}

{\em Proof.}
We assume that the fan $\Sigma$ is simplicial, which we can always
do because the formula of Proposition \ref{formula} holds for any
simplicial subdivision of faces of $\Delta^*$. 
Notice that $Ell(y,q)$ as defined in that proposition is the value  at
$\nu=0$ of 
$$
\rho(y,q,\nu)=y^{-d/2}\sum_{m\in M} \ee^{2\pi\ii  m\cdot\nu}
\left(
\sum_{n\in K^*}
y^{n\cdot {\deg} -m\cdot {\deg}^*}
q^{m\cdot n+m\cdot{\deg}^*}
G(y,q)^{d+2}
\right)
$$
where $\nu\in N_{\CC}$.
At each bidegree of $q^ky^l$ the corresponding coefficient is 
a linear combination of a finite number of exponential functions
$\ee^{2\pi\ii m\cdot\nu}$.

The idea is to split this formula for $\rho$ into the sum over 
the cones of $\Sigma$ of maximum dimension. In the smooth case it 
is precisely Bott formula.
\begin{lem}
{\rm 
For each cone $C^*\in \Sigma$ of maximum dimension we denote the
generators of its one-dimensional faces by $n_i,\,i=1,...,d+2$ suppressing
the
dependence on $C^*$. We also denote the dual basis in $M_\QQ$
by $\{m_i\}$. We choose $n_1=\deg^*$. Elements of the group $G=N/(\ZZ
n_1+...+\ZZ n_{d+2})$ could be identified with lattice points in $C^*$
whose coordinates in the basis $\{n_i\}$ are less than $1$. We call
this set of points ${\rm Box}(C^*)$. 
Then 
$$
\rho(y,q,\nu)=
\sum_{C^*\in \Sigma,\dim C^*=d+2}
{1\over|G|}
\sum_{n,l\in {\rm Box}(C^*)} y^{n\cdot \deg}
{\theta(m_1\cdot\nu,\tau)\over \theta(m_1\cdot\nu-z,\tau)}\times$$
$$
\times\prod_{i=2}^{d+2}
{\theta(-m_i\cdot\nu-m_i\cdot l-(m_i\cdot n)\tau - z)\over
\theta(-m_i\cdot\nu-m_i\cdot l-(m_i\cdot n)\tau)}.
$$
}
\label{Bott}
\end{lem}

{\em Proof of the lemma.} First of all, we must explain what the
expression above means. It is a finite sum of meromorphic functions 
of  $(\tau, z,\nu) \in H\times \CC\times N_{\CC}$. Most terms of the summation are not 
defined at $\nu=0$. We also notice that $m_1$ is orthogonal to all points
of the ${\rm Box}(C^*)$ because $\Delta^*$ is reflexive.

We start by rewriting $\rho$ as the sum over cones of any dimension
which is analogous to the formulas of Section 5. We have
$$
\rho(y,q,\nu)=y^{-{d\over2}}
\sum_{m\in M} 
\ee^{2\pi\ii m\cdot\nu}
\sum_{C^*\in \Sigma} (-1)^{{\rm codim}C^*}
\sum_{n\in C^*}
y^{n\cdot {\deg} -m\cdot {\deg}^*}
q^{m\cdot n+m\cdot{\deg}^*}
G(y,q)^{d+2}
$$
where the summation is taken over all cones in $\Sigma$ that contain
$\deg^*$. This assures that each point $n\in K^*$ contributes once.
Our next goal is to somehow get rid of cones of positive codimension.
For each of these cones there is an element $m_{C^*}\in M_1$ orthogonal to
$C^*$. Therefore, for elements of $m$ that differ by $m_{C^*}$ the
corresponding terms in the above expression for $\rho$
differ only by $\ee^{2\pi\ii m_{C^*}\cdot\nu}$. As a result,
$$\rho(y,q,\nu)
\prod_{C^*,{\rm codim} C^*>0}(1-\ee^{2\pi\ii m_{C^*}\cdot\nu})
=
y^{-{d\over2}}
\sum_{m\in M} 
\ee^{2\pi\ii m\cdot\nu}
\sum_{C^*\in \Sigma, \dim C^*=d+2} 
$$
$$
\left(
\sum_{n\in C^*}
y^{n\cdot {\deg} -m\cdot {\deg}^*}
q^{m\cdot n+m\cdot{\deg}^*}
G(y,q)^{d+2}
\right)
\prod_{C^*_1,{\rm codim} C^*_1>0}(1-\ee^{2\pi\ii m_{C^*_1}\cdot\nu})
.
$$
This identity should be understood as an identity in
$\CC[M][y,y^{-1}][[q]]$. However, we can interpret both sides 
as double series in $y$ and $q$ whose coefficients are meromorphic
functions on $N_\CC$. Really, for each $C^*$ its contribution to
$\rho$ is the supertrace over $H^0(\pi_*\msv(L),A_\CC)$ of
$y^{J[0]}q^{L[0]}\ee^{2\pi\ii A[0]}$.  It was shown in \cite{Bvertex}
that for fixed powers of $y$ and $q$ the sections of $\pi_*\msv(L)$
form a Noetherian almost-module over $\CC[C\cap M_1]$. As a result, the 
coefficients  by $y^aq^b$ in the above supertrace are Hilbert functions 
of some finitely generated modules over $\CC[C\cap M_1]$ and are therefore 
finite linear combinations of $\ee^{2\pi\ii m\cdot\nu}$ over products 
of $(1-\ee^{2\pi\ii k_im_i\cdot\nu})$ where $k_im_i$ are generators of
one-dimensional faces of $C\cap M_1$. We remark, that this happens
only after you multiply by $G(y,q)^{d+2}$, otherwise many other $m$ 
seem to contribute to a given coefficient by $y^aq^b$.

Now it remains to show that for each cone $C^*$ we have 
$$
\sum_{m\in M}
\sum_{n\in C^*}
\ee^{2\pi\ii m\cdot\nu}
y^{n\cdot {\deg} -m\cdot {\deg}^*}
q^{m\cdot n+m\cdot{\deg}^*}
G(y,q)^{d+2}=
$$
$$
={1\over|G|}
\sum_{n,l\in {\rm Box}(C^*)} y^{n\cdot \deg}
{\theta(m_1\cdot\nu,\tau)\over \theta(m_1\cdot\nu-z,\tau)}
\prod_{i=2}^{d+2}
{\theta(-m_i\cdot\nu-m_i\cdot l-(m_i\cdot n)\tau - z)\over
\theta(-m_i\cdot\nu-m_i\cdot l-(m_i\cdot n)\tau)}.
$$
There are two ways to do so. One can use the isomorphism between
${\rm BRST}_g$ cohomology of $\Fock_{M\oplus C^*}$ and the direct sum over
all
$n$ in ${\rm Box}(C^*)$ of the spaces $G$-invariant sections of the flat
space Fock space found in \cite{Bvertex}. Alternatively, a more elementary
argument is to write 
$$\sum_{m\in M}\ee^{2\pi\ii m\cdot\nu}...={1\over|G|}\sum_{l\in {\rm
Box}(C^*)}
\sum_{m\in\ZZ\{m_i\}}\ee^{2\pi\ii m\cdot\nu+m\cdot l}...$$
and then notice that 
 $$\sum_{n\in C^*}...= \sum_{n\in {\rm Box}(C^*)}\sum_{n'\in \ZZ_\geq
0\{n_i\}}...$$
Afterwards, the key formula (\ref{key}) allows us to rewrite $\rho$ 
as the sum of infinite products which are easily seen to coincide with 
ratios of theta functions. This finishes the proof of the lemma.
\hfill$\Box$

We can now go back to the proof of Proposition \ref{modular}.
It is enough to check modular properties of $Ell(y,q)$ for generators
of the Jacobi group. This follows from the transformation properties of
the theta function and Lemma \ref{Bott}. We will sketch the argument in
the hardest case of $(z,\tau)\to({z\over\tau},-{1\over\tau})$. Consider
the change in $\rho$ incurred when we change $z\to {z\over\tau}, \nu\to
{\nu\over\tau}, \tau\to -{1\over\tau}$. We have 
$$\rho({z\over\tau},-{1\over\tau},{\nu\over\tau})=
\sum_{C^*\in \Sigma,\dim C^*=d+2}
{1\over|G|}
\sum_{n,l\in {\rm Box}(C^*)} \ee^{{2\pi\ii z n\cdot \deg \over\tau}}
{\ee^{\pi\ii {(m_1\cdot\nu)^2\over\tau}}
\over
\ee^{\pi \ii {(m_1\cdot \nu-z)^2\over\tau}}
}
{\theta(m_1\cdot\nu,\tau)\over \theta(m_1\cdot\nu-z,\tau)}
\times
$$
$$
\times\prod_{i=2}^{d+2}
{\theta(-{m_i\cdot\nu\over\tau}-m_i\cdot l+{m_i\cdot n\over\tau} -
{z\over\tau},-{1\over\tau})\over
\theta(-{m_i\cdot\nu\over\tau}-m_i\cdot l+{m_i\cdot
n\over\tau},-{1\over\tau})
}=
$$
$$
=
\sum_{C^*\in \Sigma,\dim C^*=d+2}
{1\over|G|}
\sum_{n,l\in {\rm Box}(C^*)} \ee^{2\pi\ii nz\cdot \deg \over \tau}
\ee^{\pi\ii {(m_1\cdot\nu)^2\over\tau}-\pi \ii
{(m_1\cdot\nu-z)^2\over\tau}}
{\theta(m_1\cdot\nu,\tau)\over \theta(m_1\cdot\nu-z,\tau)}\times
$$
$$
\times\prod_{i=2}^{d+2}
{\theta(-m_i\cdot\nu-m_i\cdot l \tau+m_i\cdot n - z,\tau)\over
\theta(-m_i\cdot\nu-m_i\cdot l\tau+m_i\cdot n,\tau)
}
\ee^{{\pi\ii\over\tau}((m_i\cdot\nu+m_i\cdot l\tau+z-m_i\cdot
n)^2-(m_i\cdot\nu+m_i\cdot l\tau-m_i\cdot n)^2)}=
$$
$$
=
\ee^{d\pi\ii z^2\over\tau} \ee^{{2\pi\ii z\over\tau}(\deg\cdot\nu)}
\sum_{C^*}{1\over|G|}
\sum_{n,l} \ee^{2\pi\ii l\cdot \deg z} 
{\theta(m_1\cdot\nu,\tau)\over\theta(m_1\cdot\nu-z,\tau)}\times
$$
$$
\times \prod_{i=2}^{d+2}
{\theta(-m_i\cdot\nu-m_i\cdot l \tau+m_i\cdot n - z,\tau)\over
\theta(-m_i\cdot\nu-m_i\cdot l\tau+m_i\cdot n,\tau)
}.
$$
In the last step we have used that $m_1$ is orthogonal to all
elements of ${\rm Box}$ and $\sum_1^{d+2} {m_i}=\deg$.
Notice that $l$ and $n$ are switched now. Also, the change of sign
of $n$ is not important, because $n$ could be really thought of as
the element of $N/\ZZ\{n_i\}$. This implies that
$$\rho({z\over\tau},-{1\over\tau},{\nu\over\tau})=\ee^{d\pi\ii
z^2\over\tau}
\ee^{{2\pi\ii z\over\tau}(\deg\cdot\nu)}\rho(z,\tau,\nu)$$
which gives the desired modular property for $Ell(y,q)$ once we plug in
$\nu=0$. 
\hfill$\Box$

We can combine the results of Theorem \ref{modular} and Proposition
\ref{duality} to prove mirror duality of elliptic genera of Calabi-Yau
hypersurfaces in arbitrary Gorenstein toric Fano varieties.
\begin{theo}
{\rm 
Let $X$ and $X^*$ be two mirror hypersurfaces in Gorenstein toric
Fano varieties. Then 
$$Ell(X;y,q)=(-1)^d Ell(X^*;y,q).$$
}
\label{main}
\end{theo}
{\em Proof.} Use the result of \ref{duality} and transformation
property of $Ell$ under $(y,q)\to(y^{-1}q,q)$.\hfill$\Box$

\bigskip
We will now extend the results of Section 5 from smooth toric 
varieties to toric varieties with Gorenstein singularities.
We will return to the setup of that section, that is we have 
two dual lattices of $M$ and $N$ of rank $d$ and a complete fan
$\Sigma$ in $N$ which defines the toric variety $\PP$.
We first notice that (in the smooth case) the formula of
Theorem \ref{elltoric}
could be rewritten as
$$
Ell(\PP,y,q)=y^{-d/2}
\sum_{m\in M}
\sum_{C^*\in\Sigma}
(-1)^{{\rm codim} C^*}
(\sum_{n\in C^*}q^{m\cdot n}y^{\deg\cdot n})
G(y,q)^d$$
where $\deg\cdot n$ is a piecewise linear function on $N$ which 
equals $1$ on the generators of one-dimensional faces of $\Sigma$.
In general $\PP$ has only Gorenstein singularities if and only if 
this function takes integer values. This prompts the following
definition.
\begin{dfn}
{\rm
For a toric Gorenstein variety $\PP$ we define
$$
Ell(\PP,y,q)=y^{-d/2}
\sum_{m\in M}
\sum_{C^*\in\Sigma}
(-1)^{{\rm codim} C^*}
(\sum_{n\in C^*}q^{m\cdot n}y^{\deg\cdot n})
G(y,q)^d.$$
}
\label{combdef}
\end{dfn}

We remark that the expression above should be interpreted as follows.
While the series $\sum_{n\in C^*}q^{m\cdot n}y^{\deg\cdot n}$ in general
diverges at $q=0$, if we consider the product of this double series 
with $G(y,q)$ one can show that the result will only have non-negative
powers of $q$.  Alternatively, one can notice that for any
given $m$ and $C^*$ the function $\sum_{n\in
C^*}q^{m\cdot n}y^{\deg\cdot n}$ is a rational function of $q,y$ which
could be expanded around $q=0$. Besides, it could be shown that Definition
\ref{combdef} coincides with Definition \ref{EGasTrace} when one defines
$\msv(\PP)$ as in \cite{Bvertex}. The proof is very similar to the
hypersurface case. 

Since a compact toric variety is never Calabi-Yau, we do not expect
$Ell(\PP,y,q)$ to be a Jacobi form. However we will now show that
$(-1)^{d/2}Ell(\PP,-1,q)$
has expected  modular properties with respect to $\Gamma_0(2)$. The idea
is the same as in the hypersurface case. We consider
\begin{equation}
\rho(q,\nu)=
\sum_{m\in M}\ee^{m\cdot\nu}
\sum_{C^*\in\Sigma}
(-1)^{{\rm codim} C^*}
(\sum_{n\in C^*}q^{m\cdot n}(-1)^{\deg\cdot n})
G(-1,q)^d
\label{toricformula1}
\end{equation}
as a function on $H\times N_\CC$. We assume that $\Sigma$ is simplicial,
which could be done safely, because one can show that $\rho$ does not
change under crepant subdivisions of $\Sigma$. Then we rewrite it
as
\begin{equation}
\rho(q,\nu)=\sum_{C^*\in\Sigma,\dim C^*=d}{1\over|G|}
\sum_{k,l\in {\rm Box}(C^*)} (-1)^{\deg\cdot k}\prod_{i=1}^d
{\theta({1\over2}-m_i\cdot k\tau-m_i\cdot\nu-m_i\cdot l,\tau)
\over
\theta(-m_i\cdot k\tau-m_i\cdot\nu-m_i\cdot l,\tau)
}.
\label{toricformula2}
\end{equation}
Here $m_i$ denote the basis of $M_\CC$ dual to the basis of generators
$n_i$ of one-dimensional faces of $C^*$. In general $m_i\notin M$.
The proof of this formula is completely analogous to the proof of Lemma
\ref{Bott} so we skip it.

It is well-known that the group $\Gamma_0(2)$ is generated by
$\tau\to\tau+1$ and $\tau\to {\tau\over-2\tau+1}$. Clearly,
$Ell(\PP,-1,\tau)$ is not affected
by the first transformation. We will now calculate how it changes 
under the second one. 
\begin{prop}
{\rm
$$
Ell(\PP,-1,{\tau\over-2\tau+1})=(-i)^d Ell(\PP,-1,\tau)$$
}
\label{gamma02}
\end{prop}

{\em Proof.}
It is easy to derive that 
$$\theta({z\over -2\tau+1},{\tau\over -2\tau+1})=
-\ii\sqrt{2\tau-1}\,\ee^{2\pi\ii z^2\over 2\tau-1}\,\theta(z,\tau).$$
It implies 
$$
\rho({\tau\over -2\tau+1},{\nu\over -2\tau+1})
=\sum_{C^*\in\Sigma,\dim C^*=d}{1\over|G|}
\sum_{k,l\in {\rm Box}(C^*)} (-1)^{\deg\cdot k}
$$
$$
\prod_{i=1}^d
{\theta({1\over2}-m_i\cdot k{\tau\over 1-2\tau}-m_i\cdot{\nu\over
1-2\tau}-m_i\cdot l,{\tau\over 1-2\tau})
\over
{\theta(-m_i\cdot k{\tau\over 1-2\tau}-m_i\cdot{\nu\over
1-2\tau}-m_i\cdot l,{\tau\over 1-2\tau})}}
=
$$
$$
\sum_{C^*\in\Sigma,\dim C^*=d}{1\over|G|}
\sum_{k,l\in {\rm Box}(C^*)} (-1)^{\deg\cdot k}
\prod_{i=1}^d
{
\theta({-\tau+{1\over2}-m_i\cdot k\tau-m_i\cdot\nu-m_i\cdot
l(1-2\tau)\over 1-2\tau},{\tau\over 1-2\tau})  
\over
\theta({-m_i\cdot k\tau-m_i\cdot\nu-m_i\cdot
l(1-2\tau)\over 1-2\tau},{\tau\over 1-2\tau})} =
$$
$$
\sum_{C^*,\dim C^*=d}{1\over|G|}
\sum_{k,l\in {\rm Box}(C^*)} (-1)^{\deg\cdot k}
\prod_{i=1}^d
{
\theta(-\tau+{1\over2}-m_i\cdot k\tau-m_i\cdot\nu-m_i\cdot
l(1-2\tau),\tau)
\over
\theta(-m_i\cdot k\tau-m_i\cdot\nu-m_i\cdot
l(1-2\tau),\tau)}
$$
$$
\times 
\prod_{i=1}^d
\ee^{{2\pi\ii\over 2\tau-1}(-\tau+{1\over2})
(-\tau+{1\over 2}-2m_i\cdot k\tau-2m_i\cdot\nu-2m_i\cdot l(1-2\tau))}
=$$
$$
\sum_{C^*\in\Sigma,\dim C^*=d}{1\over|G|}
\sum_{k,l\in {\rm Box}(C^*)} (-1)^{\deg\cdot k}
(-\ii)^d
\prod_{i=1}^d
{
\theta({1\over2}-m_i\cdot (k-2l)\tau-m_i\cdot\nu-m_i\cdot l,\tau)
\over
\theta(-m_i\cdot (k-2l)\tau-m_i\cdot\nu-m_i\cdot l,\tau)}
$$
$$
=(-\ii)^d \rho(\tau,\nu).$$
At the last step we used the fact that we can consider $k,l$
to be representatives of the group $G$. Really, when we change 
$k$ by $n$ which is an integer combination of $n_i$, the ratio of $\theta$
functions changes by $(-1)^{m_i\cdot n}$, so overall the change
is $(-1)^{\deg\cdot n}$ which is compensated by the change in the
factor $(-1)^{\deg\cdot k}$.

Now it remains to plug in $\nu=0$ and to use (\ref{toricformula1}).
\hfill$\Box$

\begin{rem}
{\rm 
One can easily show that $Ell(\PP,-1,q)$ equals zero if the dimension
of $\PP$ is odd. Really, if we switch $k \to -k$, $l\to -l$ in the 
above summation and then use $\theta({1\over 2}-z,\tau)=-\theta(-{1\over
2}+z,\tau)=\theta({1\over 2}+z,\tau)$ we will see that
$\rho(-\nu,\tau)=(-1)^d\rho(\nu,\tau)$. So for odd $d$ 
$Ell(\PP,-1,q)=0$. When $\PP$ is smooth this could be also shown
easily by means of Chern classes.
}
\end{rem}

\begin{theo}
{\rm
Let $\PP$ be a Gorenstein toric variety of even dimension $d$.
Analogously to Section 5 we introduce
$$\widehat{Ell}(\PP,q)=(-1)^{d/2}Ell(\PP,-1,q)G(-1,q)^{-d}.$$
Then this normalized genus has transformation properties of the 
modular form of weight $d$ with respect to the group $\Gamma_0(2)$.
}
\end{theo}

{\em Proof.} We notice that 
$$G(-1,q)=\eta(2\tau)^2/\eta(\tau)^4$$
where $\eta(\tau)$ is the Dedekind $\eta$-function. Then the modular
transformation properties of $\widehat{ELL}$ follow from 
Proposition \ref{gamma02} and transformation properties of $\eta$
(cf. for example \cite{Chandra}). \hfill$\Box$

\begin{rem}
{\rm
We conjecture that $\widehat{Ell}(\PP,q)$ is a modular form. In view
of the above theorem, it simply means that it is holomorphic for
all $\tau$ and has appropriate Fourier expansions around the cusps 
of $\Gamma_0(2)$. 
}
\end{rem}


\begin{thebibliography}{99}

\bibitem{Bailey} W. N. Bailey, {\em An Algebraic Identity}, J. London
Math. Soc., {\bf 11} (1936), 156-160.

\bibitem{bat.dual} V.~V. Batyrev,  {\em Dual polyhedra and mirror
symmetry for Calabi-Yau hypersurfaces in toric varieties},
J. Algebraic Geom., {\bf 3} (1994) 493-535.

\bibitem{Bvertex} L. A. Borisov, {\em Vertex Algebras and Mirror
Symmetry}, \\ preprint math.AG/9809094.

\bibitem{Taubes} R.Bott, C.Taubes, {\em On the rigidity theorems of
Witten}, Journal of A.M.S., {\bf 2} (1989), 138-186.

\bibitem{Chandra} K. Chandrasekharan, {\em Elliptic functions}, Fundamental
Principles of Mathematical Sciences, {\bf 281}, Springer-Verlag, Berlin-New York,
1985.

\bibitem{CoxKatz} D. Cox, S. Katz, {\em  Mirror Symmetry and Algebraic
Geometry}, Mathematical Surveys and monographs, {\bf 68}, AMS, 1999.

\bibitem{Danilov} V.~I. Danilov, {\em The Geometry of Toric Varieties},
Russian Math. Surveys, {\bf 33}(1978), 97-154.

\bibitem{DMVV} R. Dijkgraaf, D. Moore, E. Verlinde, H. Verlinde, 
{\em Elliptic genera of symmetric products and second quantized strings},
Comm. Math. Phys. {\bf 185} (1997), no. 1, 197--209.

\bibitem{EOTY} T. Eguchi, H. Ooguri, A. Taormina, S.-K. Yang, 
{\em Superconformal algebras and string compactification on manifolds with 
$SU(N)$ holonomy}, Nucl. Phys. {\bf B315} (1989), 193.

\bibitem{EZ} M. Eichler, D. Zagier, {\em The theory of Jacobi forms},
Progress in Mathematics, {\bf 55}, Birkhäuser Boston, Inc., Boston, Mass.,
1985

\bibitem{Fulton} W. Fulton, {\em Introduction to toric varieties},
Princeton University Press, 1993.

\bibitem{Hirz1} F. Hirzebruch, {\em Topological methods in Algebraic Geometry},
translated from German and Appendix One by R. L. E. Schwarzenberger. With a
preface  to the third English edition by the author and Schwarzenberger. 
Appendix Two by A. Borel. Reprint of the 1978 edition. Classics in Mathematics,
Springer-Verlag, Berlin, 1995. 

\bibitem{Hirz2} F. Hirzebruch, {\em Elliptic genera of level $N$ for 
complex manifolds}, Differential Geometric methods in Theoretical Physics
(Como 1987). K. Bleuer, M. Werner Editors, NATO Adv. Sci.Inst.Ser. 
C: Math.Phys. Sci; 250. Dordrecht, Kluwer Acad.Publ.,1988.

\bibitem{Igusa} J.-I. Igusa, {\em On Siegel modular forms genus two II},
Amer. J. Math. {\bf 86} (1964), 392-412. 

\bibitem{Kac} V.~Kac, {\em Vertex algebras for beginners}, University 
Lecture Series, {\bf 10}, American Mathematical Society, Providence, RI,
1997.

\bibitem{KYY} T. Kawai, Y. Yamada, S.-K. Yang,
{\em Elliptic Genera and N=2 Superconformal Field Theory}, Nucl. Phys.
{\bf B414}(1994), 191-212.

\bibitem{Krichever} I. Krichever, {\em Generalized elliptic genera and
  Baker-Akhiezer functions}, Math. Notes, 47 (1990), 132-142.

\bibitem{Landweber.Stong} P. S. Landweber, editor, {\em Elliptic curves
and modular forms in algebraic topology}, Lecture Notes in Math.,
{\bf 1326}, Springer, Berlin, 1988.

\bibitem{LW} A. Libgober, J. Wood, {\em Uniqueness of the complex structure on
K\"{a}hler manifolds of certain homotopy types}, J. Differential Geom.
{\bf 32} (1990),
no. 1, 139--154.

\bibitem{Liu} K. Liu, {\em Modular Forms and Topology},
  Contemp. Math. {\bf 193}, AMS, 1996. 

\bibitem{MSV} F. Malikov, V. Schechtman, A. Vaintrob,
{\em Chiral de Rham complex}, preprint alg-geom/9803041.

\bibitem{Mumford} D.Mumford, {\em Tata lectures on theta. I},  with the
assistance of C. Musili, M. Nori, E. Previato and M. Stillman. Progress in
Mathematics, {\bf 28}, Birkh\"{a}user Boston, Inc., Boston, Mass., 1983.

\bibitem{Neumann}  C. D. D. Neumann, {\em The elliptic genus of Calabi-Yau $3$-
and $4$-folds, product formulae and generalized Kac-Moody algebras}, 
J. Geom. Phys., {\bf 29} (1999), no. 1-2, 5--12.

\bibitem{Oda} T. Oda, {\em Convex Bodies and Algebraic Geometry - An  
Introduction to the Theory of Toric Varieties}, Ergeb. Math. Grenzgeb.
(3), vol. 15, Springer-Verlag, Berlin, Heidelberg, New York, London,
Paris, Tokyo, 1988.

\bibitem{Totaro} B. Totaro, {\em Chern numbers of singular varieties and
  elliptic homology}, preprint, U of Chicago.

\bibitem{vanderGeer} G. van der Geer, {\em On the geometry of a Siegel
modular threefold}, Math. Ann. {\bf 260}(1982), no. 3, 317--350. 

\end{thebibliography}
\end{document}